\documentclass[11pt,sumlimits,intlimits]{amsart}
\usepackage{url,cite,amsmath,amsthm,amssymb,mathtools}
\usepackage[colorlinks=true,linkcolor=blue,citecolor=blue]{hyperref}
\usepackage[margin=1.25in]{geometry}
\usepackage{xcolor}
\usepackage{comment}
\usepackage{cleveref}
\usepackage{xypic}
\usepackage[inline]{enumitem}

\newtheorem{thmx}{Theorem}
\newtheorem{conjx}{Conjecture}

\newtheorem{theorem}{Theorem}[section]
\newtheorem*{theorem*}{Theorem}
\newtheorem{proposition}[theorem]{Proposition}
\newtheorem*{conjecture*}{Conjecture}
\newtheorem{corollary}[theorem]{Corollary}
\newtheorem{lemma}[theorem]{Lemma}
\theoremstyle{definition}
\newtheorem{definition}[theorem]{Definition}
\newtheorem{remark}[theorem]{Remark}
\newtheorem{example}[theorem]{Example}

\newcommand{\la}{\langle}
\newcommand{\ra}{\rangle}

\newcommand{\Z}{\mathbb{Z}}

\newcommand{\C}{\mathbb{C}}

\newcommand{\T}{\mathbb{T}}

\newcommand{\mK}{\mathcal{K}}

\newcommand{\Ad}{\operatorname{Ad}}

\newcommand{\id}{\operatorname{id}}

\newcommand{\dimnuc}{\operatorname{dim}_{\operatorname{nuc}}}

\newcommand{\hirsch}{\mathbf{h}}
\newcommand{\growthdegree}{\mathbf{d}}
\newcommand{\dimltc}{\dim_{\operatorname{LTC}}}

\newcommand{\fitt}{\operatorname{Fitt}}

\numberwithin{equation}{section}

\title{Nuclear dimension and virtually polycyclic groups}
\author{Caleb Eckhardt}
\address{Department of Mathematics, Miami University\newline\indent Oxford OH 45056, USA}
\email{eckharc@miamioh.edu}

\author{Jianchao Wu}
\address{Shanghai Center for Mathematical Sciences, Fudan University,\newline\indent 2005 Songhu Rd., Shanghai 200438, China}
\email{jianchao\_wu@fudan.edu.cn}

\thanks{ The first author was partially supported by an AMS-Simons research enhancement grant for PUI faculty GR000765. The second author was partially supported by NSFC Key Program No.~12231005 and National
	Key R\&D Program of China 2022YFA100700. }

\makeatletter
\@namedef{subjclassname@2020}{
	\textup{2020} Mathematics Subject Classification}
\makeatother


\begin{document}
\maketitle
\begin{abstract} We show that the nuclear dimension of a (twisted) group C*-algebra of a virtually polycyclic group is finite. This prompts us to make a conjecture relating finite nuclear dimension of group C*-algebras and finite Hirsch length, which we then verify for a class of elementary amenable groups beyond the virtually polycyclic case. In particular, we give the first examples of finitely generated, non-residually finite groups with finite nuclear dimension. A parallel conjecture on finite decomposition rank is also formulated and an analogous result is obtained. Our method relies heavily on recent work of Hirshberg and the second named author on actions of virtually nilpotent groups on $C_0(X)$-algebras. 

\end{abstract}
\tableofcontents

\section{Introduction}
Nuclear dimension plays an essential role in the Elliott classification program for simple nuclear C*-algebras. Indeed finite nuclear dimension is a necessary condition for a C*-algebra to be classifiable in the sense of \cite{Gong20,Carrion23}. On the other hand Winter and Zacharias--the architects of nuclear dimension--also showed connections of nuclear dimension to dynamical systems and coarse geometry  and declared their hope to see its role in dimension type conditions in other areas of noncommutative geometry \cite{WinterZacharias2010nuclear}.

The classification program was certainly one motivation for this work, but we are at least as interested in a group theoretic characterization of those groups whose C*-algebras have finite nuclear dimension.
The Hirsch length of a group (see \Cref{sec:Polycyclicgroups}),  plays the role of dimension for a vector space (for example it is increasing with respect to inclusion and obeys a rank-nullity condition). A group with finite Hirsch length--as it was originally defined--is called \emph{virtually polycyclic}.   
Our main result shows a connection between finite Hirsch length and finite nuclear dimension.  Indeed  our nuclear dimension bound  is in terms of Hirsch length.

\begin{thmx}[{see Theorems~\ref{thm:finitend} and~\ref{thm:twistednucdim}}] \label{thm:A} Let $G$ be a virtually polycyclic group and $\sigma$ a circle-valued 2-cocycle on $G$.  Then the twisted group C*-algebra $C^*(G,\sigma)$ has finite nuclear dimension. In the case of untwisted group C*-algebras of such groups, there is a finite bound on their nuclear dimension that depends only on the Hirsch length.
\end{thmx}

The definition of Hirsch length was relaxed by Hillman \cite{Hillman91} to an extent that many elementary amenable, non virtually polycyclic groups have finite Hirsch length while retaining its meaning as a dimension theory. Hillman's Hirsch length agrees with the classic definition for virtually polycyclic groups. Examples of non virtually polycyclic groups with finite nuclear dimension (e.g. the lamplighter group) are found in \cite{Hirshberg17}  and all of these examples have finite Hirsch length.  
On the other hand the only known examples of elementary amenable groups with infinite nuclear dimension all have infinite Hirsch length (see \Cref{sec:wr-infinite}).  We feel \Cref{thm:A} lends strong evidence for the following 
\begin{conjx} \label{conj:dimnuc} A finitely generated elementary amenable group has finite Hirsch length if and only if its group C*-algebra has finite nuclear dimension.
\end{conjx}

To give further evidences to this conjecture, 
we notice that all previously known examples of finitely generated groups with finite nuclear dimension happen to be residually finite.  We use the techniques developed here to fill this small gap in the literature and to showcase how our methods extend beyond the virtually polycyclic world. 

\begin{thmx}[{see \Cref{thm:nonRFexample}}] \label{thm:B}  There are infinitely many finitely generated solvable non-residually finite groups whose group C*-algebras have finite nuclear dimension.  All of these examples have finite Hirsch length as defined by Hillman. 
\end{thmx}

Furthermore, 
the following result confirms \Cref{conj:dimnuc} for a class of elementary amenable groups that includes not only all virtually polycyclic groups but also many examples with infinite Hirsch length. 

\begin{thmx} \label{thm:C} 
	Let $K$ be a virtually polycyclic group and let $H$ be a finitely generated virtually nilpotent group. 
	Then the following are equivalent: 
	\begin{enumerate}
		\item \label{thm:C::nilpotent} The wreath product $K \wr H$ has finite Hirsch length.   
		\item \label{thm:C::finite} Either $H$ or $K$ is finite.  
		\item \label{thm:C::dimnuc} The C*-algebra $C^*(K \wr H)$ has finite nuclear dimension.
	\end{enumerate}
\end{thmx}

Note that the class of wreath product groups in \Cref{thm:C} contains many non-residually finite and non-virtually solvable examples. Indeed, if $H$ is infinite and $K$ contains a nonabelian simple subgroup $S$, then the presence of the infinite direct sum $S^{\oplus H}$ inside $K \wr H$ prevents the latter from being virtually solvable, while it also follows from \cite[Theorem~3.1]{Gruenberg1957} that $K \wr H$ is not residually finite. 

Finally, before we outline our proof methods, we record a few words about decomposition rank versus nuclear dimension. 
For any C*-algebra $A$, the \emph{decomposition rank} (see \cite{Kirchberg04}) of $A$ is trivially greater than or equal to $\dimnuc(A).$ Hence  any decomposition rank bound immediately yields a nuclear dimension bound. If $G$ is a finitely generated virtually nilpotent group, then $C^*(G)$ has finite decomposition rank \cite{EckhardtGillaspyMcKenney2019Finite} and therefore finite nuclear dimension.

In this work we focus mainly on {nuclear dimension} because many polycyclic groups have infinite decomposition rank. It is possible that every finitely generated non virtually nilpotent group has infinite decomposition rank.  In fact there are special cases where we have complete information. Let  $\alpha$ be an automorphism of $\Z^d.$   Then $\dimnuc(C^*(\Z^d\rtimes_\alpha \Z))<\infty$ by \cite{Hirshberg17} and the following are equivalent: 
\begin{enumerate}
\item \label{before-conj:II::nilpotent} The semidirect product $\Z^d\rtimes_\alpha \Z$ is virtually nilpotent.
\item \label{before-conj:II::eigenvalue} Every eigenvalue of $\alpha$ has modulus 1.
\item \label{before-conj:II::sqd} The C*-algebra $C^*(\Z^d\rtimes_\alpha \Z)$ is strongly quasidiagonal. 
\item \label{before-conj:II::dr} The C*-algebra $C^*(\Z^d\rtimes_\alpha \Z)$ has finite decomposition rank.
\end{enumerate}
The equivalence of \eqref{before-conj:II::nilpotent} and \eqref{before-conj:II::eigenvalue} is classic.
The implications \eqref{before-conj:II::nilpotent} $\Rightarrow$ \eqref{before-conj:II::dr} $\Rightarrow$ \eqref{before-conj:II::sqd} $\Rightarrow$ \eqref{before-conj:II::nilpotent} are respectively  \cite[Theorem 6.6]{EckhardtGillaspyMcKenney2019Finite}, \cite[Theorem 5.3]{Kirchberg04} and \cite[Theorem 3.3]{Eckhardt15}. In other words  the groups $\Z^d\rtimes_\alpha\Z$ all have finite nuclear dimension but have finite decomposition rank precisely when they are virtually nilpotent. 
This led us to propose the following conjecture in parallel to \Cref{conj:dimnuc}: 

\begin{conjx} \label{conj:dr} A finitely generated elementary amenable group is virtually nilpotent if and only if its group C*-algebra has finite decomposition rank.
\end{conjx}

Parallel to \Cref{thm:C} above, we also have the following observation that verifies \Cref{conj:dr} for a class of wreath products. 

\begin{thmx} \label{thm:D} 
	Let $K$ and $H$ be finitely generated virtually nilpotent groups. Then the following are equivalent: 
	\begin{enumerate}
		\item \label{thm:D::nilpotent} The wreath product $K \wr H$ is virtually nilpotent.   
		\item \label{thm:D::finite} Either $H$ is finite or $K$ is trivial. 
		\item \label{thm:D::sqd} The C*-algebra $C^*(K \wr H)$ is strongly quasidiagonal.  
		\item \label{thm:D::dr} The C*-algebra $C^*(K \wr H)$ has finite decomposition rank.
	\end{enumerate}
\end{thmx}

\begin{proof}
	The implications \eqref{thm:D::nilpotent} $\Rightarrow$ \eqref{thm:D::dr} $\Rightarrow$ \eqref{thm:D::sqd} $\Rightarrow$ \eqref{thm:D::finite} are respectively \cite[Theorem 6.6]{EckhardtGillaspyMcKenney2019Finite}, \cite[Theorem 5.3]{Kirchberg04} and \cite[Theorem 3.4]{CarrionDadarlatEckhardt13}. 
	The implication \eqref{thm:D::finite} $\Rightarrow$ \eqref{thm:D::nilpotent} is immediate from the standard permanence properties of virtually nilpotent groups. 
\end{proof}

\subsection{Outline of Proof}
Our proof of \Cref{thm:A}  draws on the methods of \cite{EckhardtMcKenney2018Finitely,EckhardtGillaspyMcKenney2019Finite} that show finite decomposition rank of virtually nilpotent groups and incorporates the recent results of \cite{HirshbergWuActions} on nuclear dimension of virtually nilpotent actions on $C_0(X)$-algebras.

Let $G$ be a finitely generated nilpotent group.  Then $C^*(G)$ decomposes as a continuous field over the Pontryagin dual of its center.  The broad tactic of \cite{EckhardtMcKenney2018Finitely} was to obtain a uniform bound on the fibers and then to apply an estimate of Carri\'{o}n \cite{Carrion11} to bound the decomposition rank of $C^*(G)$.

If $G$ is virtually polycyclic then it may have no (useful) interpretation as a continuous field. Hence the methods of \cite{EckhardtMcKenney2018Finitely} are no longer immediately useful.
But $G$ does contain a nilpotent normal subgroup  $N\trianglelefteq G$, called the \emph{Fitting subgroup}, such that $G/N$ is virtually abelian. So we may decompose $C^*(G)$ as a twisted crossed product of $C^*(N)$ by $G/N$. This suggested we could use the recent work of  Hirshberg and the second named author in \cite{HirshbergWuActions} to attack the virtually polycyclic case. To do this, however, we need to generalize the main result in \cite{HirshbergWuActions} to the case of twisted crossed products as follows (see \Cref{thm:HW} for the precise statement):

Let $X$ be a locally compact Hausdorff space with finite covering dimension. Suppose $G$ is finitely generated, virtually nilpotent and admits a twisted action on a $C_0(X)$ algebra $A$ such that $G$ restricts to an action on $C_0(X)$. For each $x\in X$, let $G_x$ be the stabilizer.  If there is a uniform bound on the nuclear dimension of twisted crossed products $A_x\rtimes H$ as $x$ ranges through $X$ and $H$ ranges through subgroups of $G_x$, then
we can show that 
the twisted crossed product $A\rtimes G$ has finite nuclear dimension.

Coming back to the problem of bounding the nuclear dimension of $C^*(G) \cong C^*(N) \rtimes G/N$, we take $X$ to be the dual of the center of $N$. 
Essentially \Cref{thm:HW} reduces a complicated dynamical problem to the more tractable problem of dealing with actions of the stabilizers on the fibers.  
 
More precisely, it remains to find a uniform bound on the nuclear dimension of the twisted crossed products $C^*(N)_{\chi} \rtimes H/N$, where $\chi$ ranges over all characters of the center of $N$ and $H/N \leq (G/N)_{\chi}$. We observe that $C^*(N)_{\chi} \rtimes H/N$ is itself a quotient of $C^*(H)$. 
At this point, we need to discuss two possible cases: the first case is when this C*-quotient factors through a quotient group of $H$ with a lower Hirsch length, and the second is when this fails. 
Although the second case is more complicated, we end up getting an optimal result, namely, said C*-quotient is a direct sum of unital simple separable C*-algebras with nuclear dimension at most $1$. 
To this end, \Cref{lem:strongouterfitting} is the key lemma connecting group theory to operator algebras as it shows how the maximality of the Fitting subgroup forces certain actions to be strongly outer, after which we get to apply, as in \cite{EckhardtMcKenney2018Finitely}, the work of Matui and Sato on strongly outer actions \cite{MatuiSato2014}. 
With the second case resolved, the first case naturally leads us to adopt an inductive argument on the Hirsch length. 
Thus the nuclear dimension bound we obtain will necessarily be uniform with regard to Hirsch length. 

Similar but simpler applications of the main theorem in \cite{HirshbergWuActions} yield the finite nuclear dimension results in \Cref{thm:B} and \Cref{thm:C}. 

The converse direction of \Cref{thm:C} amounts to showing infinite nuclear dimension for wreath products of (certain) infinite groups. Here we need to adopt a different strategy than the usual one for groups like $\mathbb{Z}^m \wr \mathbb{Z}^n$ that relies on the residual finiteness of $\mathbb{Z}^n$, in order to circumvent a potential problem resulting from the dimension reduction phenomenon of homogeneous C*-algebras with UHF fibers \cite{TikuisisWinter2014}. A more detailed discussion of this difficulty can be found in \Cref{example:ZwrZ}. Our strategy is centered around keeping track of the various copies of the augmentation ideal $I(H) \mathop{\triangleleft} C^*(H)$ in an infinite tensor product $C^*(H)^{\otimes \infty}$ (see \Cref{lemma:embedding-ideal}). 

\subsection{Organization of the Paper} After reviewing basic facts about nuclear dimension, virtually nilpotent groups, Hirsch length, and, more substantially, the theory of twisted group C*-algebras and twisted crossed products in \Cref{sec:prelim}, 
we establish in \Cref{sec:extension-twisted} two results regarding the relationship between central extensions and (twisted) group C*-algebras, to be used in later proofs. 
Then we discuss in \Cref{sec:mathcalC} a class of nice C*-algebras (namely, unital, separable, simple, with unique trace, and having finite nuclear dimension) that is well-suited for a bootstrap argument useful in our nuclear dimension bounds. In \Cref{sec:FND}, we give finite nuclear dimension bounds in the various settings stated in our main theorems (\ref{thm:A}, \ref{thm:B}, and \ref{thm:C}). As mentioned above, our proofs here rely on a recent result of Hirshberg and the second named author \cite{HirshbergWuActions}, which we upgrade to the twisted setting in \Cref{thm:HW}. Finally in \Cref{sec:wr-infinite}, we prove a large class of infinite-by-infinite wreath product groups give rise to group C*-algebras with infinite nuclear dimension, thus completing the proof of \Cref{thm:C}. 

\subsection{Acknowledgments} We are grateful to the anonymous referee for a number of helpful suggestions. In particular, one suggestion led to a more conceptual treatment of \Cref{lem:centralextension}.

\section{Preliminaries} \label{sec:prelim}
\noindent We collect preliminaries about nuclear dimension, groups, cohomology and twisted crossed products.

\subsection{Nuclear dimension} \label{sec:dimnuc}
We recall the definition and basic properties of nuclear dimension. 

\begin{definition}
	A completely positive map $\varphi \colon A \to B$ is \emph{order zero} if $\varphi(a) \varphi(a') = 0$ for any $a, a' \in A_+$ with $a a' = 0$.

	The \emph{nuclear dimension} of a $C^*$-algebra $A$ is the infimum of all the natural numbers $d$ such that for any finite subset $F \subseteq A$ and any $\varepsilon > 0$ there are a finite dimensional $C^*$-algebra $B$, a completely positive contractive map $\psi \colon A \to B$, and completely positive order zero maps $\varphi^{(0)}, \varphi^{(1)}, \ldots, \varphi^{(d)}  \colon B \to A$, such that $\| \sum_{l=0}^d \varphi^{(l)} \circ \psi (a) - a \| < \varepsilon$ for any $a \in F$.  
\end{definition}

Throughout the paper, we write $\otimes$ for the spatial (as opposed to maximal) tensor product of C*-algebras. 

\begin{lemma}[{see \cite[Proposition~2.3, Corollary~2.8(i), Proposition~2.9, and Corollary~2.10 ]{WinterZacharias2010nuclear} and \cite[Lemma~3.1]{Carrion11}\footnote{See also \cite[the discussion above Lemma~3.3]{Hirshberg17}}}] \label{lemma:dimnuc-basic}
	Let $A$, $B$, and $C$ be C*-algebras, and let $(A_i)_{i \in I}$ be a net of C*-algebras. Then the following hold: 
	\begin{enumerate}
		\item \label{lemma:dimnuc-basic::continuous-trace} If $A$ is a separable continuous trace algebra, then $\dimnuc(A)$ agrees with the covering dimension of the spectrum of $A$. 
		\item \label{lemma:dimnuc-basic::extension} If $0 \to A \to B \to C \to 0$ is a short exact sequence of C*-algebras, then we have $\max \{\dimnuc (A) , \dimnuc (C)\} \leq \dimnuc (B) \leq \dimnuc (A) + \dimnuc (C) + 1$. 
		\item \label{lemma:dimnuc-basic::limit} We have $\dimnuc (\varinjlim_{i \in I} A_i) \leq \liminf_{i \in I} \dimnuc (A_i)$. 
		\item \label{lemma:dimnuc-basic::stabilization} We have $\dimnuc(A \otimes \mathcal{K}) = \dimnuc(A)$, where $\mathcal{K}$ is the algebra of compact operators on a separable Hilbert space. 
		\item \label{lemma:dimnuc-basic::field} If $A$ is a separable $X$-C*-algebra for some locally compact metric space $X$, then we have $\dimnuc(A) + 1 \leq (\dim(X) + 1) \cdot \sup_{x \in X} (\dimnuc(A_x)+1)$.  
	\end{enumerate}
\end{lemma}

\subsection{Polycyclic groups} \label{sec:Polycyclicgroups}
We refer the reader to Dan Segal's book \cite{Segal1983Polycyclic} for more information on virtually polycyclic groups.

 A group $G$ is called \textbf{polycyclic} if it has a normal series
\begin{equation}
G_0=\{ e \}\trianglelefteq G_1\trianglelefteq\cdots \trianglelefteq G_n\trianglelefteq G_{n+1}=G \label{eq:normalseries}
\end{equation}
such that each quotient group $G_{i+1}/G_i$ is cyclic. The \textbf{Hirsch length} of a polycyclic group, defined as $\hirsch(G)=|\{ i\in\{ 0,...,n \}:G_{i+1}/G_i \textup{ is infinite} \}|,$ is independent of the normal series \eqref{eq:normalseries}.  If $G$ is virtually polycyclic, then one defines $\hirsch(G)=\hirsch(N)$ where $N$ is any finite index polycyclic subgroup of $G$. 
In particular, we have $\hirsch(\mathbb{Z}^n \oplus F) = n$ if $F$ is finite. 
It follows that for a finitely generated abelian group $G$, if we write $\widehat{G}$ for its Pontryagin dual, we have 
\begin{equation} \label{eq:Hirsch-Pontryagin}
	\hirsch(G) = \dim(\widehat{G})\; . 
\end{equation}

The class of (virtually) polycyclic groups is closed under subgroups, quotients and extensions. 
One notices that 
all finitely generated nilpotent groups are polycyclic by considering either the upper or lower central series, 
and all polycyclic groups are finitely generated. 

In \cite[Theorem~1]{Hillman91}, the notion of Hirsch length is extended to work for all elementary amenable groups, 
so that if $G$ is elementary amenable and $N$ is a normal subgroup, then 

\begin{equation} \label{eq:Hirsch-length-additive}
	\hirsch(G)=\hirsch(N)+\hirsch(G/N) \; .
\end{equation}

The \textbf{Fitting subgroup} of a group $G$ is the group generated by all normal, nilpotent subgroups of $G$ and is denoted by $\fitt(G)=\la N\trianglelefteq G: N \textup{ is nilpotent}  \ra.$  If $G$ is virtually polycyclic, then $\fitt(G)$ is nilpotent \cite[Corollary 1.B.13]{Segal1983Polycyclic} and $G/\fitt(G)$ is virtually abelian.  Key to our investigations is the fact that every virtually polycyclic group $G$ contains a characteristic, torsion free, finite index subgroup $K$ such that $K/\fitt(K)$ is free abelian \cite[Page~92]{Segal1983Polycyclic}. 
In particular, by taking the upper or lower central series of $\fitt(K)$ in the above, we see that $G$ has a finite chain of characteristic subgroups $1 = G_n \trianglelefteq \ldots \trianglelefteq G_1 \trianglelefteq G_0 \trianglelefteq G$ such that $G / G_0$ is finite and $G_{i-1} / G_i$ is finitely generated abelian for $i = 1, 2, \ldots, n$. 

The reason we want to impose torsion-freeness here is related to our need to arrange strongly outer actions (see \Cref{def:inner-outer} and \Cref{lem:strongouterfitting}).

Let $N$ be a group.  Its upper central series is defined inductively as $Z_0(N) = \{ e \}$ and $Z_{i+1}(N)/Z_i(N) = Z(N/Z_i(N)).$  If $N$ is nilpotent then there is an integer $c$ so $N=Z_c(N).$  Moreover for any $x\in Z_i(N)$ with $i\geq1$ and $y\in N$ we have $[x,y] \in Z_{i-1}(N).$   When $N$ is torsion free and nilpotent then $N/Z_i(N)$ is torsion free for all $i \geq 0$ (see \cite[Theorem~1.2]{Jennings55}).  By a slight modification of the proof of \cite[Theorem~1.2]{Jennings55} we obtain the following fact needed for the proof of \Cref{lem:strongouterfitting}. 

\begin{lemma} \label{lem:uniqueroot} Let $N\trianglelefteq K$ such that $N$ is a torsion free nilpotent normal subgroup of $K.$ Let $x\in K, s\in N$ and $n\in \Z\setminus \{ 0 \}$ such that $[x,s^n]=1.$  Then $[x,s]=1.$
\end{lemma}
\begin{proof} First notice that $[x,s^n] = 1$ if and only if $[x,s^{-n}]=1$ so we may suppose that $n\geq 1.$

For any elements $a,b,c$ in a group we have the identity $[a,bc] = [a,c][c,[b,a]][a,b].$ With $a=x, b=s$ and $c=s^{n-1}$ we have
\begin{equation} \label{eq:commutatorpower}
[x,s^n] = [x,s^{n-1}][s^{n-1},[s,x]][x,s]
\end{equation}
Suppose that $[x,s] \neq 1.$ Since $N$ is normal in $K$ we have $[x,s] \in N.$ By way of contradiction let $i\geq 0$ be the largest non-negative integer such that $[x,s]\not\in Z_i(N).$ Since $[s,x] = [x,s]^{-1} \in Z_{i+1}(N)$ we have 
$[s^{k},[s,x]] \in Z_i(N)$ for any $k \in \mathbb{Z}$. Hence by induction we obtain from (\ref{eq:commutatorpower}) that $[x,s^n]Z_i(N) = [x,s]^nZ_i(N).$  Since $[x,s^n] = 1$ and $N/Z_i(N)$ is torsion free we have $[x,s] \in Z_i(N)$ a contradiction.
\end{proof}

\subsection{Group cohomology} \label{sec:groupcohomology}
In this section we recall the ideas from group (co)homology necessary for an understanding of the proof of \Cref{thm:twistednucdim}.  Our presentation is extremely brief, very specific to our needs, unmotivated and self-contained. We refer the reader to $\S$~II.3 and $\S$~III.1 of \cite{Brown1982Cohomology} for the full treatment.  

 Let $G$ be a group. For each $n\geq1$ we define $C_n(G)$ to be the free abelian group with basis  $\{ (g_1,...,g_n):g_i\in G \}$. The \textbf{boundary operators} $\partial_n:C_n(G)\rightarrow C_{n-1}(G)$ for $n=2,3$ are the homomorphisms defined by
\begin{equation} \label{eq:boundary3}
\partial_3(g_1,g_2,g_3)=(g_2,g_3)-(g_1g_2,g_3)+(g_1,g_2g_3)-(g_1,g_2) \; ,
\end{equation} 
\begin{equation}
\partial_2(g_1,g_2)=(g_1)-(g_1g_2)+(g_2) \; .
\end{equation}
Elements of $Z_2(G):=\textup{ker}(\partial_2)$ are called \textbf{2-cycles}.  Elements of $B_2(G):=\textup{Im}(\partial_3)$ are called \textbf{2-boundaries}.  
For any $g, h \in G$, applying $\partial_3$ to $(g,e,h)$ and $(g, g^{-1}, g)$ in $C_3(G)$, we get 
\begin{equation} \label{eq:2-cycles-examples}
	(g, e) - (e, h) , \  (g , g^{-1}) - (g^{-1} , g) \in B_2(G)  .
\end{equation}
One easily calculates that $\partial_2\partial_3$ is the zero map, whence $B_2(G)\leq Z_2(G).$  The \textbf{second homology group} of $G$ is defined as $H_2(G)=Z_2(G)/B_2(G)$. 
One also defines the \textbf{first homology group} $H_1(G) = C_1(G) / B_1(G)$, where $B_1(G):=\textup{Im}(\partial_2)$. 
A useful fact here is that if $G$ is finite, then $H_1(G)$ and $H_2(G)$ are torsion groups. 

For each $n\geq 1$, one defines $C^n(G,\T):=\textup{Hom}(C_n(G),\T)$ and adopts multiplicative notations. The \textbf{coboundary operators} $\partial^n:C^{n-1}(G,\T)\rightarrow C^n(G,\T)$ for $n=2,3$ are defined by
\begin{equation*}
\partial^n(u)(x)=u(\partial_n(x))
\end{equation*}
Elements of $Z^2(G,\T):=\textup{ker}(\partial^3)$ are called \textbf{2-cocycles} (or more precisely, \emph{circle-valued 2-cocycles}).  Elements of $B^2(G,\T):=\textup{Im}(\partial^2)$ are called \textbf{2-coboundaries}.  One easily calculates that $\partial^2\partial^3$ is the zero map, whence $B^2(G,\T)\leq Z^2(G,\T).$  The \textbf{second cohomology group of} $G$ \textbf{with coefficients in }$\T$ is defined as $H^2(G,\T)=Z^2(G,\T)/B^2(G,\T).$ We say that two cocycles $\omega,\omega'$ are \textbf{cohomologous} if $\omega B^2(G,\T)=\omega'B^2(G,\T).$

Let $\omega\in Z^2(G,\T)\leq \textup{Hom}(C_2(G),\T)$ be a 2-cocycle.  Then $\omega$ is completely determined by the values $\omega(g_1,g_2)$ for $g_1,g_2\in G.$  Hence one may consider a 2-cocycle $\omega$ as a function $\omega:G\times G\rightarrow \T$ and the condition $\omega\in \textup{ker}(\partial^3)$ provides the so-called \textbf{cocycle identity}
\begin{equation}\label{eq:cocycleidentity}
\omega(g_1,g_2g_3)\omega(g_2,g_3)=\omega(g_1g_2,g_3)\omega(g_1,g_2)\quad \textup{ for all }g_1,g_2,g_3\in G
\end{equation}

On the other hand, for any $\omega\in \textup{Hom}(C_2(G),\T)$, it is clear that $\omega\in Z^2(G,\T)$ if and only if $B_2(G)\leq \textup{ker}(\omega)$ if and only if $\omega$ factors through a group homomorphism from $C_2(G) / B_2(G)$ to $\T$.
This gives rise to a canonical isomorphism 
\begin{equation}\label{eq:Z2_as_hom}
	Z^2(G,\T) \xrightarrow{\cong} \textup{Hom}(C_2(G) / B_2(G),\T)
\end{equation}
Moreover, when we view a 2-cocycle $\omega$ as a homomorphism $C_2(G) / B_2(G) \to \T$ and restrict its domain to $H_2(G) = Z_2(G) / B_2(G) \leq C_2(G) / B_2(G)$, we see that the restricted homomorphism only depends on the cohomology class $[\omega] \in H^2(G,\T)$. 
This gives rise to a canonical homomorphism 
\begin{equation}\label{eq:H2_as_hom}
	H^2(G,\T) \xrightarrow{\cong} \textup{Hom}(H_2(G),\T)
\end{equation}
which, thanks to the divisibility of $\T$, is also bijective. 
Indeed, this follows from the more general UCT sequence $0 \to \mathrm{Ext}(H_1(G), M) \to H^2(G, M) \to \textup{Hom}(H_2(G),M) \to 0$ for arbitrary abelian groups $M$ and the fact that $\mathrm{Ext}(H_1(G), M) = 0$ when $M$ is divisible. It can also be seen directly: injectivity follows from the fact that any homomorphism $B_1(G) \to \mathbb{T}$ extends to $C_1(G)$, since this means that any homomorphism $C_2(G) / B_2(G) \to \mathbb{T}$ that vanishes on $Z_2(G) / B_2(G)$ is in the image of $\partial^1$; surjectivity follows from the fact that any homomorphism $Z_2(G) / B_2(G) \to \mathbb{T}$ extends to $C_2(G) / B_2(G)$ (while injectivity relies on the divisibility of $\T$, surjectivity does not, since $C_2(G) / B_2(G)$ splits as a direct sum $Z_2(G) / B_2(G) \oplus B_1(G)$, a fact we exploit in the proof of \Cref{lem:centralextension}).

\subsection{Twisted crossed products} \label{sec:twisted-croe}

Suppose $N \trianglelefteq G$ are groups. If the sequence $0\rightarrow N\rightarrow G\rightarrow G/N\rightarrow 0$ splits, then one may decompose $C^*(G)$ as a crossed product of $C^*(N)$ by $G/N.$  If the short exact sequence does not split,  then the twisted crossed products of \cite{PackerRaeburn1989twisted} provide the technical framework to study $C^*(G)$ in terms of the ``simpler" groups $N$ and $G/N$. We refer the reader to Packer and Raeburn's \cite{PackerRaeburn1989twisted} for more information on the general case or to \cite{Bedos91} for a brisker introduction to the discrete case.

Let $A$ be a C*-algebra and $G$ a discrete group. Let $UM(A)$ denote the unitary group of the multiplier C*-algebra of $A$.  For any $u\in UM(A)$ we define the automorphism $\Ad u$ by $\Ad u(x)=uxu^*$. A \textbf{twisted dynamical system} $(G,A,\alpha,\omega)$  are maps $\alpha:G\rightarrow \textup{Aut}(A)$ and $\omega:G\times G\rightarrow UM(A)$ satisfying for all $s,t,r\in G$, 
 \begin{enumerate}
 \item[$\cdot$]  $\omega(s,e)=\omega(e,s)=1$ and $\alpha_e=\id$; 
 \item[$\cdot$] $\alpha_s\circ \alpha_t=\Ad(\omega(s,t))\circ\alpha_{st}$;
 \item[$\cdot$] $\alpha_r(\omega(s,t))\omega(r,st)=\omega(r,s)\omega(rs,t)$. 
 \end{enumerate}
Note that the last equation is a more general version of the cocycle identity \eqref{eq:cocycleidentity}. A genuine action $\alpha \colon G \curvearrowright A$ thus corresponds to a twisted dynamical system $(G,A,\alpha,\omega)$ where $\omega(s,t) = 1$ for any $s,t \in G$. 

 A \textbf{covariant representation} $(\pi,U)$ of $(G,A,\alpha,\omega)$ is defined as a representation $\pi:A\rightarrow B(\mathcal{H})$ and a map $U:G\rightarrow U(B(\mathcal{H}))$ satisfying
 \begin{equation*}
 \pi\circ\alpha_s=\Ad U_s\circ\pi, \textup{ and }U_sU_t=\pi(\omega(s,t))U_{st}
 \end{equation*}

 When $A$ is unital, the \textbf{algebraic twisted crossed product} $A\rtimes_{\alpha,\omega, \mathrm{alg}}G$ is the $*$-algebra generated by $A$ and unitaries $(u_g)_{g\in G}$ subject to the relations
 \begin{equation*}
 	u_gu_h=\omega(g,h)u_{gh}, \textup{ and } u_g a u_g^* =\alpha_g(a) \; .
 \end{equation*}
 When $A$ is non-unital, $A\rtimes_{\alpha,\omega, \mathrm{alg}}G$ is the ideal generated by $A$ in $A^+\rtimes_{\alpha,\omega, \mathrm{alg}}G$, that is, it is spanned by the products $a u_s$, where $a \in A$ and $s \in G$. 
 The \textbf{maximal twisted crossed product} $A\rtimes_{\alpha,\omega, \mathrm{max}}G$ is the C*-envelope of $A\rtimes_{\alpha,\omega, \mathrm{alg}}G$. 
 Note that there is a canonical way to identify covariant representations of $(G,A,\alpha,\omega)$, $*$-representations of $A\rtimes_{\alpha,\omega, \mathrm{alg}}G$, and $*$-representations of $A\rtimes_{\alpha,\omega, \mathrm{alg}}G$. 

Let $\pi:A\rightarrow B(\mathcal{H})$ be a faithful representation of $A.$  The regular representation $(\tilde{\pi},u)$ of $(G,A,\alpha,\omega)$ is the covariant representation on $\ell^2(G,\mathcal{H})$ defined by
\begin{equation*}
\tilde{\pi}(a)(\xi)(t)=\pi(\alpha_t(a))\xi(t), (u_s\xi)(t)=\pi(\omega(s,t))\xi(ts) \; .
\end{equation*}
The \textbf{reduced twisted crossed product} $A\rtimes_{\alpha,\omega, \mathrm{r}}G$ is
the C*-algebra generated by $\tilde{\pi}(A)$ and $\{ u_t:t\in G \}$. The definition is independent of the faithful representation $\pi$. 
When $G$ is amenable, the canonical $*$-homomorphism $A\rtimes_{\alpha,\omega, \mathrm{max}}G \to A\rtimes_{\alpha,\omega, \mathrm{r}}G$ is a $*$-isomorphism \cite[Theorem~3.11]{PackerRaeburn1989twisted}, so in this case we shall denote both by $A\rtimes_{\alpha,\omega}G$ for the sake of brevity.

Two twisted dynamical systems $(G,A,\alpha,\omega)$ and $(G,A,\beta,\sigma)$ are said to be \emph{exterior equivalent} if there are $v_g \in UM(A)$ such that for any $s, t \in G$ we have 
 \begin{enumerate}
 	\item[$\cdot$] $\beta_s = \Ad v_s \circ \alpha_s$, and 
 	\item[$\cdot$] $\sigma(s,t) = v_s \, \alpha_s(v_t) \, \omega(s,t) \, v_{st}^*$. 
 \end{enumerate}
This induces a $*$-isomorphism 
\begin{align*}
A \rtimes_{\alpha,\omega, \mathrm{alg}}G \xrightarrow{\cong} A \rtimes_{\beta, \sigma, \mathrm{alg}}G 
\, , \quad  a \, u_g \mapsto a \, v_g \, u_g 
\; ,
\end{align*}
which extends to $*$-isomorphisms 
\[
A \rtimes_{\alpha,\omega, \mathrm{max}}G \xrightarrow{\cong} A \rtimes_{\beta, \sigma, \mathrm{max}}G 
\quad \text{ and } \quad 
A \rtimes_{\alpha,\omega, \mathrm{r}}G \xrightarrow{\cong} A \rtimes_{\beta, \sigma, \mathrm{r}}G 
\; .
\]

\subsection{The stabilization trick}\label{sec:stabilization}
In \cite{PackerRaeburn1989twisted}, Packer and Raeburn also developed a technique that turns the stabilization of a twisted crossed product into an untwisted crossed product in an explicit way. This provides a convenient way to generalize results about untwisted crossed product to twisted ones. We recall it below.  

Let $(G,A,\alpha,\omega)$ be a twisted dynamical system  
and let $\mK(\ell^2(G))$ denote the C*-algebra of compact operators on $\ell^2(G)$. 
Recall that $A \otimes \mK(\ell^2(G))$ is naturally isomorphic to the C*-algebra $\mK(A \otimes \ell^2(G))$ of compact operators on the (right) Hilbert $A$-module $A \otimes \ell^2(G)$, while its multiplier C*-algebra is naturally isomorphic to the C*-algebra $B(A \otimes \ell^2(G))$ of adjointable $A$-module maps. 
In particular, $M(A)$ canonically embeds into $M(A \otimes \mK(\ell^2(G)))$ by acting on the first tensor factor of $A \otimes \ell^2(G)$, while any $*$-automorphism $\varphi$ of $A$ extends canonically to $M(A \otimes \mK(\ell^2(G)))$ via taking conjugation by the map $\varphi \otimes \id_{\ell^2(G)}$, so that $\varphi (a \otimes T) = \varphi(a) \otimes T$ for any $a \in A$ and $T \in \mK(\ell^2(G))$.  

For any $s \in G$, we let $\lambda^\omega_s \in B(A \otimes \ell^2(G))$ be the unitary defined by 
\[
	\lambda^\omega_s (a \otimes \delta_g) = {\omega(s, g)}^* a \otimes \delta_{sg} \quad \text{ for any } a \in A \text{ and } g \in G \; . 
\] 
Then the cocycle identity for $(G,A,\alpha,\omega)$ implies  
\[
	\lambda^\omega_s \, \alpha_s  (\lambda^\omega_t)  = \lambda^\omega_{s t} \, \omega(s,t)^* \quad \text{ for any } s,t \in G \; .
\]
It follows that the automorphisms 
\[
	\alpha_s^\omega = \Ad( \lambda^\omega_s ) \circ \alpha_s 
	\in  \textup{Aut} (M (A \otimes \mK(\ell^2(G)))) 
	\quad \text{ for } s \in G
\]
restricts to a genuine action $\alpha^\omega$ of $G$ on $A \otimes \mK(\ell^2(G))$ that is exterior equivalent to $(G, A  \otimes \mK(\ell^2(G)), \alpha, \omega \otimes 1)$. 
Hence we have 
$*$-isomorphisms 
\[
	(A\rtimes_{\alpha,\omega, \mathrm{max}}G) \otimes \mK(\ell^2(G)) \cong (A \otimes \mK(\ell^2(G)) ) \rtimes_{\alpha,\omega \otimes 1, \mathrm{max}}G \xrightarrow{\cong} (A \otimes \mK(\ell^2(G)) )\rtimes_{\alpha^\omega, \mathrm{max}}G 
\]
and 
\[
	(A\rtimes_{\alpha,\omega, \mathrm{r}}G) \otimes \mK(\ell^2(G)) \cong (A \otimes \mK(\ell^2(G)) ) \rtimes_{\alpha,\omega \otimes 1, \mathrm{r}}G \xrightarrow{\cong} (A \otimes \mK(\ell^2(G)) )\rtimes_{\alpha^\omega, \mathrm{r}}G 
	\; .
\]

An immediate application of these $*$-isomorphisms is the fact that if $G$ is amenable, then the canonical $*$-epimorphism from the maximal twisted crossed product to the reduced one is a $*$-isomorphism. 

\subsection{Twisted $G$-$X$-C*-algebras}\label{sec:G-X-algebras}
Our techniques require us to discuss how C*-algebras decompose over (subalgebras of) their centers in ways compatible with (twisted) actions. 

Let $(G,A,\alpha,\omega)$ be a twisted dynamical system with $G$ discrete. Note that the extension of $\alpha$ to the multiplier algebra $M(A)$ restricts to a genuine action on the center $ZM(A)$. Now let $X$ be a locally compact Hausdorff space. 
A twisted $G$-$X$-C*-algebra consists of a twisted dynamical system $(G,A,\alpha,\omega)$, a (genuine) action $\beta \colon G \curvearrowright C_0(X)$ (sometimes still denoted by $\alpha$; this is equivalent to having an action of $G$ on $X$ by homeomorphisms), and a $G$-equivariant homomorphism $\iota \colon C_0(X) \to ZM(A)$ that is nondegenerate in the sense that $\overline{\iota (C_0(X)) A} = A$. 

Given a twisted $G$-$X$-C*-algebra $(G,A,\alpha,\omega, X, \beta, \iota)$ and any $x \in X$, we write $A_x$ for the quotient C*-algebra $A / (\overline{\iota (C_0(X \smallsetminus \{x\})) A})$, sometimes called the \textbf{fiber} at $x$, and, for any $a \in A$, write $a_x$ for its image under the quotient map. If $x$ is fixed by a subgroup $H$ of $G$, then the ideal $\overline{\iota (C_0(X \smallsetminus \{x\})) A}$ is preserved by the automorphisms $\alpha_h$ for $h \in H$, and thus we have a twisted dynamical system $(H, A_x, \alpha^x_{|H}, \omega^x_{|H})$, where $\alpha^x_{|H}$ is induced from the $(\alpha_h)_{h \in H}$ and $\omega^x_{|H}(s,t) = \omega(s,t)_x$ for any $s, t \in H$. 
Now for any twisted dynamical system $(G, A, \alpha', \omega')$ that is exterior equivalent to $(G,A,\alpha,\omega)$ (say, via $(v_s)_{s \in G}$), since clearly $\alpha$ and $\alpha'$ induce the same action on $ZM(A)$, it follows that $(G,A,\alpha',\omega', X, \beta, \iota)$ is again a twisted $G$-$X$-C*-algebra. Moreover, for any $x \in X$ that is fixed by a subgroup $H$ of $G$, the induced twisted dynamical systems $(H, A_x, \alpha^x_{|H}, \omega^x_{|H})$ and $(H, A_x, \alpha'{}^x_{|H}, \omega'{}^x_{|H})$ are clearly also exterior equivalent (via $((v_s)_x)_{s \in G}$).

An important special case here is when the space $X$ is a homogeneous space $G/H$ for some subgroup $H \leq G$, with the action $\beta \colon G \curvearrowright C_0(G/H)$ given by left translation. In this case, we have a form of Green's imprimitivity isomorphism (\cite{Green1980imprimitivity}): 
\begin{equation} \label{eq:imprimitivity}
	\begin{array}{rcll}
		A \rtimes_{\alpha, \omega, \mathrm{r}} G &\cong& \left(A_{[e]} \rtimes_{\alpha, \omega, \mathrm{r}} H\right) \otimes \mathcal{K} (\ell^2(G/H)) 
		&\text{ and } \\
		A \rtimes_{\alpha, \omega, \mathrm{max}} G &\cong& \left(A_{[e]} \rtimes_{\alpha, \omega, \mathrm{max}} H\right) \otimes \mathcal{K} (\ell^2(G/H))  & .
	\end{array}
\end{equation}
More precisely, upon choosing a set-theoretic lift $c \colon G/H \to G$ and identifying $A$ with the direct sum $\bigoplus_{[g] \in G/H} A_{[g]}$, 
we simply map any $(u_h a) \otimes e_{[g], [g']} \in \left(A_{[e]} \rtimes_{\alpha, \omega, \mathrm{max}} H\right) \otimes \mathcal{K} (\ell^2(G/H))$ to $u_{c([g])} \, u_h a \, u_{c([g'])}^* \in A \rtimes_{\alpha, \omega, \mathrm{max}} G$, and conversely, map any $u_{g} a \in A \rtimes_{\alpha, \omega, \mathrm{max}} G$, where $a \in A_{[g']}$ (and thus $u_{g} a = u_{c([gg'])} \,   u_{c([gg'])^{-1} g \, c([g'])} \, v \,  \alpha_{c([g'])^{-1}}(a) \,u_{c([g'])}^*$ where $v \in UM(A)$ satisfies $u_{c([gg'])^{-1} g \, c([g'])} \, v = u_{c([gg'])}^* \, u_g \, u_{c([g'])}$), to $\left( u_{c([gg'])^{-1} g \, c([g'])} \, v \, \alpha_{c([g'])^{-1}}(a) \right) \otimes e_{[gg'], [g']} \in \left(A_{[e]} \rtimes_{\alpha, \omega, \mathrm{max}} H\right) \otimes \mathcal{K} (\ell^2(G/H))$. 
These assignments will generate a pair of mutually inverse isomorphisms between the maximal completions. 
This pair of isomorphisms also descend to the reduced level, since upon fixing a faithful representation $A_{[e]} \to B(\mathcal{H})$, one can represent both sides, via regular representations of the crossed products, on the same Hilbert space $\mathcal{H} \otimes \ell^2(G/H) \otimes \ell^2(H) \otimes \ell^2(G/H)$, such that the above isomorphisms intertwine the two representations.

\subsection{Twisted group C*-algebras}\label{sec:twistedgroupCalgebras}
Twisted group C*-algebras are a special case of twisted crossed products, but it is sometimes clunky to describe them as such, so they receive their own section.

Let $G$ be a discrete group and $\omega \in Z^2(G,\T)$.  
Define operators $u_g\in B(\ell^2(G))$ by
\begin{equation*}
u_g(\delta_h) = \omega(g,h)\delta_{gh} \text{ for }h\in G \; .
\end{equation*}  
Let $\C[G,\omega] = \text{span}\{ u_g:g\in G \}.$ We then have
$u_gu_h=\omega(g,h)u_{gh}$, from which we obtain 
\[ 
	u_e = \omega(e,e) \cdot 1 \text{ and } u_g^* = \overline{\omega(g,g^{-1}) \omega(e,e)} u_{g^{-1}} \; .
\]
The \textbf{maximal twisted group C*-algebra} $C^*(G,\omega)$ is the $C^*$-envelope of $\C[G,\omega]$, while
the \textbf{reduced twisted group C*-algebra} $C_{\mathrm{r}}^*(G,\omega)$ is the norm closure of $\C[G,\omega]$ in $B(\ell^2(G))$. The \textbf{twisted von Neumann algebra} $L(G,\omega)$ is the von Neumann algebra generated by $\C[G,\omega]$ in $B(\ell^2(G))$.  The map $\tau(x) = \la x\delta_e,\delta_e \ra$ defines a faithful trace on $L(G,\omega)$. 
Note that if two 2-cocycles $\omega$ and $\omega'$ are cohomologous, which is equivalent to saying that there is a function $f \colon G \to \T$ such that 
\[
	f(g g') \omega (g, g') = f(g) f(g')  \omega' (g, g') \quad \text{ for any } g , g' \in G \; ,
\] 
then there is a $*$-isomorphism 
\[
\C[G,\omega] \xrightarrow{\cong} \C[G,\omega'] \, , \quad u_g \mapsto f(g) u_g \; ,
\]
which extends to $*$-isomorphisms $C_{\mathrm{r}}^*(G,\omega) \cong C_{\mathrm{r}}^*(G,\omega')$ and $L(G,\omega) \cong L(G,\omega')$.  

A 2-cocycle is \textbf{normalized} if $\omega(g,g^{-1}) = 1$ for all $g\in G.$ It is well-known (see for example \cite{Kleppner62}) that any 2-cocycle is cohomologous to a normalized cocycle so we may assume all of our cocycles are normalized and in particular assume that

\[ 
	u_e = 1 \text{ and } u_g^* =  u_{g^{-1}} \; .
\]

\subsection{Twisted group C*-algebras as twisted crossed products}\label{sec:groupCalgebrastwistedcp} 
We will find it convenient to decompose twisted group C*-algebras as twisted crossed products with respect to normal subgroups and vice versa. We recall this common procedure in this section. For the sake of simplicity, we focus on the maximal completions, though a version for reduced completions also exists. 

Let $G$ be a discrete group, $\sigma\in Z^2(G,\T)$ and $N\trianglelefteq G$ a normal subgroup.  
Let $c \colon G/N\rightarrow G$ be a lifting of the quotient homomorphism such that $c(e_{G/N})=e_G.$ 
We denote the restriction of $\sigma$ to $N \times N$ still by $\sigma$ 
and define $\omega:G/N\times G/N\rightarrow U(C^*(N,\sigma))$ and $\alpha:G/N\rightarrow \textup{Aut}(C^*(N,\sigma))$ as
 \begin{equation} \label{eq:tgastcp-action}
 \omega(s,t)=u_{c(s)}u_{c(t)}u_{c(st)}^{*}, \textup{ and } \alpha_s=\text{Ad}(u_{c(s)}).
 \end{equation}
 Then $(G/N, C^*(N,\sigma),\alpha,\omega)$ defines a twisted dynamical system. The covariant representation $(\id \colon C^*(N,\sigma)\rightarrow C^*(N,\sigma)$, $t\mapsto u_{c(t)})$ provides the isomorphism
 \begin{equation} \label{eq:tgastcp}
 C^*(N,\sigma)\rtimes_{\alpha,\omega,\mathrm{max}}G/N \cong C^*(G,\sigma) \; .
 \end{equation}
 
More generally, if $(G, A, \beta, \sigma)$ is a twisted dynamical system, then using the same formulas as in \eqref{eq:tgastcp-action}, we also have an isomorphism
\begin{equation} \label{eq:tcpastcp}
	(A \rtimes_{\beta, \sigma,\mathrm{max}} N)\rtimes_{\alpha,\omega,\mathrm{max}}G/N \cong A \rtimes_{\beta, \sigma,\mathrm{max}} G \; .
\end{equation}

\begin{example} \label{ex:central-extension}
	As a special case of the construction above, consider a central extension $1 \to N \to G \to G/N \to 1$. Then combining the decomposition above and the Gelfand duality, we obtain a $*$-isomorphism 
	\[
	C^*(G) \cong C(\widehat{N}) \rtimes_{\id, \sigma,\mathrm{max}} (G/N) \; ,
	\]
	where $\widehat{N}$ is 
	the Pontryagin dual of $N$, $\id$ stands for the trivial automorphisms, and $\sigma \colon G/N \times G/N \to U(C(\widehat{N}))$ satisfies 
	\[
	\sigma ([s], [t]) (\chi) = \chi \left(c([s])c([t])c([st])^{-1} \right) \quad \text{ for any } s,t \in G \text{ and } \chi \in  \widehat{N} \; ,
	\]
	where $c \colon G/N \rightarrow G$ is some identity preserving lift. 
	It follows that as an $\widehat{N}$-C*-algebra (in the obvious way), 
	$C^*(G)$ satisfies that 
	the fiber at an arbitrary $\chi \in \widehat{N}$ is isomorphic to $C^*(G/N, \sigma_\chi)$, where $\sigma_\chi$ is defined as $\sigma(-,-)(\chi)$. 
	We shall see in \Cref{lem:centralextension} that in fact, every twisted group C*-algebra arises this way.   
\end{example}

The following generalization of \Cref{ex:central-extension} plays an essential role in our proof of \Cref{thm:finitend}. 

\begin{example} \label{ex:beyond-central-extension}
	Consider two normal subgroups $M$ and $N$ in a group $G$ such that $M \leq Z(N)$. 
	Following \Cref{ex:central-extension} and choosing some identity preserving lift $d \colon N/M \to N$,  
	we view $C^*(N)$ as an $\widehat{M}$-C*-algebra, such that for each $\chi\in\widehat{M}$, 
	we have 
	\[
		C^*(N)_\chi \cong C^*(N/M,\sigma_\chi) \quad \text{ with } \sigma_\chi ([s], [t]) = \chi \left(d([s])d([t])d([st])^{-1} \right)  \text{ for any } s,t \in N  \; .
	\]
	Then, choosing another identity preserving lift $c \colon G/N \to G$, we consider the twisted action $(\alpha,\omega)$ of $G/N$ on $C^*(N)$ as defined in \eqref{eq:tgastcp-action} such that
	\begin{equation*}
	C^*(G) \cong C^*(N) \rtimes_{\alpha,\omega,\mathrm{max}} G/N \; .
	\end{equation*}
	This twisted action provides a (right) action $\beta$ of $G/N$ by homeomorphisms on $\widehat{M}$ which is dual to the conjugation action of $G$ on ${M}$ factored through $N$ (clearly $N$ acts trivially), that is, we have 
	\[
		\beta_{[g]} (\chi) (s) = \chi (g s g^{-1})  \quad \text{ for any } g \in G , \text{ any } \chi \in \widehat{M} ,  \text{ and any } s \in M \; .
	\]
	Viewing $\beta$ as a (left) genuine action on $C(\widehat{M})$, we have $\alpha_{[g]} (a) = \beta_{[g]} (a)$ for any $g \in G$ and $a \in C^*(M) \cong C(\widehat{M})$. 
	This turns $C^*(N)$ into a twisted $G/N$-$\widehat{M}$-C*-algebra in the sense of \Cref{sec:G-X-algebras}, whose twisted crossed product recovers $C^*(G)$. 
	
\end{example}

\section{Twisted group C*-algebras and central extensions} \label{sec:extension-twisted}
In this section, we prepare for later use two results (\Cref{lem:centralextension} and \Cref{prop:central-by-finite}) regarding central extensions of groups. 
The first result will be used in \Cref{thm:twistednucdim} and, along with its preliminary version (\Cref{lem:centralextension-canonical}), can be understood as a converse statement to \Cref{ex:central-extension}: 
while \Cref{ex:central-extension} decomposes the group C*-algebra of any central extension of a group $G$ into twisted group C*-algebras of $G$, \Cref{lem:centralextension} shows that every twisted group C*-algebra of $G$ arises this way.
This more-or-less follows from the proof of \cite[Corollary 1.3]{PackerRaeburn1992structure}. However, it is simpler to provide a straightforward proof 
rather than explain how it follows from \cite[Corollary 1.3]{PackerRaeburn1992structure}, hence we do that.

In the following we canonically identify $Z^2(G,\T)$ with $\widehat{C_2(G)/B_2(G)}$ as in \eqref{eq:Z2_as_hom}, and $H^2(G,\T)$ with $\widehat{H_2(G)}$ as in \eqref{eq:H2_as_hom}. 

\begin{lemma} \label{lem:centralextension-canonical} 
	Let $G$ be a discrete group. Then there is a canonical central extension 
	\begin{equation*}
		1 \rightarrow C_2(G)/B_2(G) \rightarrow E_G \rightarrow G\rightarrow 1
	\end{equation*}
	such that for any $\omega\in Z^2(G,\T) \cong \widehat{C_2(G)/B_2(G)}$, the twisted group C*-algebra $C^*(G,\omega)$ is canonically isomorphic to $C^*(E_G)_{\omega}$, that is, the quotient of $C^*(E_G)$ obtained by viewing $C^*(E_G)$ as an $\widehat{C_2(G)/B_2(G)}$-C*-algebra as in \Cref{ex:central-extension}. 
\end{lemma}

\begin{proof}
	We define $E_G$ to be the product set $(C_2(G)/B_2(G)) \times G$ with a binary operation defined by 
	\begin{equation} \label{lem:centralextension-canonical::eq:multiplication}
		(x_1,g_1)\cdot (x_2,g_2)=( x_1+x_2+ [(g_1,g_2)] - [(e,e)] ,g_1g_2) \quad \text{ for any } (x_1,g_1), (x_2,g_2) \in E_G
	\end{equation}
	By \eqref{eq:boundary3} we see that for all $g_1,g_2,g_3\in G$ we have $(g_2,g_3)-(g_1g_2,g_3)+(g_1,g_2g_3)-(g_1,g_2) \in B_2(G)$. From this fact it follows that the operation is associative.  
	It follows from \eqref{eq:2-cycles-examples} that 
	$( 0 , e)$ serves as the identity of said operation and 
	$(x,g)^{-1} = (-x - [(g, g^{-1})] + [(e,e)], g^{-1})$ for any $(x,g) \in E_G$.
	Therefore $E_G$ is a group under the above operation. 
	
	One checks that the embedding $C_2(G)/B_2(G) \to  E_G$ taking $x$ to $\mathring{x} := (x , e)$ is a homomorphism, whose image is $(C_2(G)/B_2(G)) \times\{ e \}$. 
	Moreover, it follows from \eqref{eq:2-cycles-examples} that 
	\begin{equation} \label{lem:centralextension::eq:central-multiplication}
		\mathring{y} \cdot (x, g) = (x + y , g) = (x ,g ) \cdot \mathring{y} \quad \text{ for any } x, y \in C_2(G)/B_2(G) \text{ and } g \in G \; ,
	\end{equation}	
	whence $(C_2(G)/B_2(G)) \times\{ e \}$ is contained in the center of $E_G$. It is clear that $G \cong E_G / ((C_2(G)/B_2(G)) \times\{ e \})$. This verifies that we have a central extension $1 \rightarrow C_2(G)/B_2(G) \rightarrow E_G \rightarrow G \rightarrow 1$. 
	
	Now fix $\omega\in Z^2(G,\T)$. 
	In order to identify the quotient $C^*(E_G)_\omega$ with $C^*(G,\omega)$, we apply \Cref{ex:central-extension}. 
	To this end, we let $c \colon G \to E_G$ be the identity preserving lift taking $g$ to $(0,g)$. A direct computation yields, for any $g, h \in G$, 
	\begin{equation*}\label{lem:centralextension::eq:lift}
		c(g) c(h) c(gh)^{-1} = ([(g,h)] - [(e,e)], e) \ \text{ and thus } \omega \left(c(g) c(h) c(gh)^{-1}\right) = \omega(g,h) \overline{\omega(e,e)} \; ,
	\end{equation*}
	the latter equation being equivalent to 
	\[
		\omega \left( (z c(g)) (z c(h)) (z c(gh))^{-1}\right) = \omega(g,h) \quad \text{ where } z = \omega(e,e) \; .
	\] 
	It follows from \Cref{ex:central-extension} that there is a canonical $*$-isomorphism $C^*(G, \omega) \to C^*(E_G)_\omega$ taking each generator $u_g \in C^*(G, \omega)$ to $\omega(e,e) [u_{c(g)}] \in C^*(E_G)_\omega$. 
\end{proof}

\begin{remark}
	As suggested by the anonymous referee, there is an alternative but equivalent construction of the above canonical central extension that is well adapted to operator algebraic applications. Given any $\omega\in Z^2(G,\T)$, the twisted product group $\T \times_\omega G$ consists of pairs $(z, g) \in \T \times G$, with multiplication defined by $(z, g) (w, h) = (z w \, \omega(g,h), gh)$ (note that the identity is given by $(\overline{\omega(e,e)}, e)$). Then the group $E_G$ is canonically isomorphic to the subgroup in $\prod_{\omega\in Z^2(G,\T)} \T \times_\omega G$ generated by $\left\{ (1, g)_{\omega \in  Z^2(G,\T)} \colon g \in G \right\}$. 
	
	To see this, one considers the canonical embedding $E_G \hookrightarrow \prod_{\omega\in Z^2(G,\T)} \T \times_\omega G$ that maps $(x,g) \in (C_2(G)/B_2(G)) \times G$ as in the proof of \Cref{lem:centralextension-canonical} to $(\omega(x) \overline{\omega(e,e)}, g)_{\omega \in  Z^2(G,\T)} \in \prod_{\omega\in Z^2(G,\T)} \T \times_\omega G$. The image of this embedding agrees with the subgroup described above, because $E_G$ is generated by $\left\{ ([(e,e)], g) \colon g \in G \right\}$, a fact that one can prove simply by noting that an arbitrary generator $[(g,h)]$ of $C_2(G)/B_2(G)$ can be written inside $E_G$ as $([(e,e)], g) ([(e,e)], h) ([(e,e)], gh)^{-1}$. 
\end{remark}

Next we shall refine the previous result by replacing the kernel $C_2(G)/B_2(G)$ with a typically much smaller group, namely $H_2(G)$. This will be important for our nuclear dimension bound in \Cref{thm:twistednucdim}. A minor price we pay here is the canonicality of the construction: we shall see in \Cref{appendix:celh} that the central extension below is not always canonical, but we can precisely quantify how non-canonical it is (see \Cref{remark:H2-extension-uniqueness} and \Cref{remark:H2-extension-uniqueness-weak}; while we find this discussion interesting and potentially useful to future endeavors, our main results do not depend on it so we relegated this material to the appendix). 

\begin{proposition} \label{lem:centralextension} Let $G$ be a discrete group. Then there is a (possibly non-canonical) central extension 
	\begin{equation*}
		1 \rightarrow H_2(G)\rightarrow E\rightarrow G\rightarrow 1
	\end{equation*}
	such that for any $\omega\in Z^2(G,\T)$, the twisted group C*-algebra $C^*(G,\omega)$ is a quotient of $C^*(E)$. In fact, if we view $C^*(E)$ as an $\widehat{H_2(G)}$-C*-algebra as in \Cref{ex:central-extension} and view $[\omega]$ as an element in $\widehat{H_2(G)} = \textup{Hom}(H_2(G),\T)$ as in \eqref{eq:H2_as_hom}, then we have 
	\[
	C^*(G,\omega) \cong C^*(E)_{[\omega]} \; .
	\]
\end{proposition}
\begin{proof} 
	We shall construct the desired central sequence as a quotient of the one in \Cref{lem:centralextension-canonical}, so that they fit into a commutative diagram 
	\begin{equation}\label{lem:centralextension:eq:central_sequences}
		\xymatrix{
			1 \ar[r] & C_2(G)/B_2(G) \ar[r] \ar@{->>}[d]^\pi & E_G \ar[r] \ar@{->>}[d]^\pi & G \ar[r] \ar@{=}[d] & 1 \\
			1 \ar[r] & H_2(G) \ar[r] & E \ar[r] & G \ar[r] & 1 
		}
	\end{equation}
	where the vertical maps are surjective. 
	To this end, consider the exact sequence
	\begin{equation*}
		0\rightarrow H_2(G)  \xrightarrow{\iota} C_2(G)/B_2(G) 
		\xrightarrow{\underline{\partial_2}} B_1(G)\rightarrow0.
	\end{equation*}
	which is obtained as a quotient of the exact sequence $0\rightarrow Z_2(G)\rightarrow C_2(G)\xrightarrow{\partial_2} B_1(G)\rightarrow0$. 
	Since $B_1(G)$ is free abelian, the above exact sequence splits. 
	More precisely, upon choosing a lifting homomorphism $\sigma \colon B_1(G)\rightarrow C_2(G) / B_2(G)$ (i.e., it satisfies $\underline{\partial_2} \circ \sigma=\id$), we obtain a direct sum decomposition $ H_2(G) \oplus B_1(G) \cong  C_2(G) / B_2 $ implemented by $\iota + \sigma$ in one direction and $\pi \oplus \underline{\partial_2}$ in the other, where $\pi =\id - \sigma \circ \underline{\partial_2} \colon C_2(G) / B_2(G) \to H_2(G)$. 
	
	Since $\sigma(B_1(G))$ is a central subgroup of $E_G$, we may form the quotient group $E := E_G / \sigma(B_1(G))$. 
	Since $ \sigma(B_1(G)) = \ker(\pi) $, the associated quotient map $E_G \to E$ extends $\pi \colon C_2(G) / B_2(G) \to H_2(G)$, whence we use the same notation $\pi$. This gives rise to the desired commutative diagram involving central extensions in \eqref{lem:centralextension:eq:central_sequences}. 
	
	Now fix $\omega\in Z^2(G,\T) \cong \widehat{C_2(G)/B_2(G)}$. We define $\omega' = \omega \circ \pi \in \widehat{C_2(G)/B_2(G)}$. Since $\pi \circ \iota = \mathrm{id}$ on $H_2(G)$, it follows that $\omega$ and $\omega'$ restricts to the same homomorphism on $H_2(G)$, whence $[\omega] = [\omega']$ in $H^2(G, \T)$ by \eqref{eq:H2_as_hom}. It follows that $C^*(E)_{[\omega]} = C^*(E)_{[\omega']}$, and $C^*(G,\omega) \cong C^*(G,\omega')$ by \Cref{sec:twistedgroupCalgebras}. On the other hand, it follows from \Cref{ex:central-extension} that the quotient map $C^*(E_G) \to C^*(G)$ induced by $\pi \colon E_G \to E$ in turn induces a $*$-isomorphism $C^*(E_G)_{\omega'} \cong C^*(E)_{[\omega']}$. Combining these isomorphisms shows $C^*(G,\omega) \cong C^*(E)_{[\omega]}$. 
\end{proof}

Now we turn to the second result, which roughly says that when we apply \Cref{ex:central-extension} to the case of central $\mathbb{Z}^d$-extensions of finite groups, the resulting twisted C*-algebra fibers can be untwisted in a locally trivial way. 
This will be used in \Cref{thm:dimnuc-wr-infty}.  

\begin{proposition} \label{prop:central-by-finite}
	Let $G$ be a discrete group and let $M$ be a finite-index central subgroup in $G$ that is isomorphic to $\mathbb{Z}^d$. 
	Then $C^*(G)$ is $*$-isomorphic to the C*-algebra of continuous sections on a locally trivial $C^*(G/M)$-bundle. 
\end{proposition}

\begin{proof}
	Let us write $F := G/ M$ and fix a set-theoretic lifting $c \colon F \to G$ as in \Cref{ex:central-extension}. 
	The assignment $(s,t) \mapsto c(s)c(t)c(st)^{-1} \in M$, for any $s, t \in F$, induces a homomorphism $\varphi \colon H_2(F) \to M$. 
	Since $F$ is finite, $H_2(F)$ is torsion. As $M$ is torsion-free, we see that 
	$\varphi = 0$.

	By \Cref{ex:central-extension}, we see that $C^*(G) \cong C(\widehat{M}) \rtimes_{\id, \sigma, \mathrm{max}} F$ and it is an $\widehat{M}$-C*-algebra with the fiber at each $\chi \in \widehat{M}$ isomorphic to $C^*(F, \sigma_\chi)$, where $\sigma_\chi$ is defined by
	\[
		\sigma_\chi (s, t) = \chi \left( c(s)c(t)c(st)^{-1} \right) \quad \text{ for any } s,t \in F \text{ and } \chi \in  \widehat{M} \; .
	\]
	Note that the 2-cocycle $\sigma_\chi$ induces $[\sigma_\chi] \in \widehat{H_2(F)}$ that agrees with $\chi \circ \varphi$. It then follows from \eqref{eq:H2_as_hom} that $\sigma_\chi$ is cohomologous to $0$, i.e., $\sigma_\chi \in B^2(F, \mathbb{T})$, and thus also $C^*(F, \sigma_\chi) \cong C^*(F)$. 
	
	It remains to show local triviality, that is, any $\chi \in \widehat{M}$ has an open neighborhood $W$ such that $C_0(W) \cdot C^*(G) \cong C_0(W) \otimes C^*(F)$ as $W$-C*-algebras. 
	To this end, recall that by definition, $B^2(F, \mathbb{T})$ is the image of $\partial^2 \colon C^1(F, \mathbb{T}) \to C^2(F, \mathbb{T})$. Under the canonical identifications $C^1(F, \mathbb{T}) \cong \mathbb{T}^F$ and $C^2(F, \mathbb{T}) \cong \mathbb{T}^{F^2}$, both groups become compact Lie groups and $\partial^2$ is a continuous group homomorphism. It follows that $B^2(F, \mathbb{T})$ is also a compact Lie group and the subspace topology on it coming from $C^2(F, \mathbb{T})$ agrees with the quotient topology coming from $C^1(F, \mathbb{T})$. 
	It is a classical result (see \cite{Mostow}) that quotients of compact Lie groups have local cross-sections; in this case, it means that for
	any $\omega \in B^2(F, \mathbb{T})$
	there is an open neighborhood $V_\omega \subseteq B^2(F, \mathbb{T})$ 
	and a continuous map $\lambda_\omega \colon V_\omega \to C^1(F, \mathbb{T})$ such that $\partial^2 \circ \lambda_\omega = \id_{V_\omega}$. 
	 
	Now for any $\chi \in \widehat{M}$, we let $W$ be the pre-image of $V_{\sigma_\chi}$ under the continuous map $\sigma \colon \widehat{M} \to B^2(F, \mathbb{T})$, $\chi' \mapsto \sigma_{\chi'}$, and define a continuous map $\beta := \lambda_{\sigma_\chi} \circ \sigma|_W \colon W \to C^1(F, \mathbb{T})$. 
	Then for any $\chi' \in W$ and any $s, t \in F$, we have $\sigma_{\chi'} (s,t) = \left( \partial^2 \circ \beta (\chi') \right) (s,t) = \beta (\chi') (s) \, \beta (\chi') (t) \, \overline{\beta (\chi') (st)}$.  
	It follows that $C_0(W) \cdot C^*(G) \cong C_0(W) \otimes C^*(F)$ via the following pair of mutually inverse $*$-isomorphisms: 
	\begin{align*}
		\Phi \colon &C_0(W) \cdot C^*(G) \to C_0(W) \otimes C^*(F)  && \text{and} & \Psi \colon C_0(W) \otimes C^*(F) &\to C_0(W) \cdot C^*(G)\\
		&f \cdot u_g \mapsto \left( f \cdot u_{g \left( c([g]) \right)^{-1}} \cdot {\beta (-) ([g])}\right) \otimes u_{[g]} &&& f \otimes u_{s} & \mapsto f \cdot \overline{\beta (-) (s)} \cdot u_{c(s)}
	\end{align*}
	where $f \in C_0(W)$, $g \in G$, and $s \in F$, and we note that $u_{g \left( c([g]) \right)^{-1}} \in C^*(M)$, the latter being identified with $C(\widehat{M})$. 
	This completes the proof. 
\end{proof}

\begin{example}
	Note that the fiber bundle in \Cref{prop:central-by-finite} is typically not globally trivial. For example, if $G = \mathbb{Z}$ and $M = 2 \mathbb{Z}$, then $C^*(G) \cong C(\mathbb{T})$ is the C*-algebra of continuous sections on the M\"{o}bius bundle $\mathbb{T} \xrightarrow{ \times 2} \mathbb{T}$ over $\widehat{M} \cong \mathbb{T}$. 
\end{example}

 \section{A useful class of C*-algebras} \label{sec:mathcalC}
In this section, we isolate and prove some theorems about a nice class of C*-algebras.  This allows the proof of the main theorem to flow smoothly.
\begin{definition}\label{def:mathcalC}
Let $\mathcal{C}$ be the class of all unital, separable, simple  C*-algebras with a unique tracial state and nuclear dimension less than or equal to 1.
\end{definition}
\begin{remark}\label{rem:equivalentCdef}  Let $A$ be a unital, separable, simple, infinite-dimensional and nuclear C*-algebra with  a unique tracial state.  Combining results of \cite{Winter2012nuclear} and \cite{Bosa2019Covering}, the following are equivalent
\begin{enumerate}
\item $\dimnuc(A)\leq 1$
\item $\dimnuc(A)<\infty$
\item $A\otimes \mathcal{Z}\cong A$ where $\mathcal{Z}$ is the Jiang-Su algebra. 
\end{enumerate}
Hence in the definition of $\mathcal{C}$ we may replace nuclear dimension less than or equal to 1 with either finite nuclear dimension or, in the case when $A$ is infinite-dimensional, $\mathcal{Z}$-stability.

Also note that if we write $\tau$ for the unique trace of $A$, then the von Neumann algebra $M := \pi_\tau(A)''$ is an injective tracial factor (\cite{EffrosLance1977}) and thus isomorphic to either $M_n$ for some positive integer $n$ or the hyperfinite II$_1$ factor $\mathcal{R}$ (\cite{Connes1976}). 
\end{remark}
We recall examples of elements of $\mathcal{C}$ and permanence properties of the class that will be used to prove our main result.
\begin{theorem}\textup{(\cite{EckhardtGillaspyMcKenney2019Finite})} \label{thm:nilpotenttwist}
Let $N$ be a finitely generated nilpotent group and $\pi$ an irreducible representation of $N.$  Then the C*-algebra generated by $\pi(N)$ is an element of $\mathcal{C}.$
In particular, if $N$ is a torsion free nilpotent group and $\chi\in\widehat{Z(N)}$ is a faithful character, then
 $C^*(N/Z(N),\sigma_\chi)\in \mathcal{C}.$ 
\end{theorem}
\begin{remark} There is no hope to extend \Cref{thm:nilpotenttwist} to the virtually polycyclic case.
Let $G$ be a virtually polycyclic group and $J$ a maximal ideal of $C^*(G).$  Then $C^*(G)/J$ satisfies all of the conditions of being in $\mathcal{C}$ except possibly having unique trace. We give a specific example. Consider $G = \Z^2 \rtimes \Z$ where the action is implemented by $A = \left[ \begin{array}{cc} 2 & 1\\ 1 & 1  \end{array} \right].$    The group C*-algebra $C^*(\Z^2\rtimes \Z)$ has a simple $\mathcal{Z}$-stable quotient that admits more than one trace.

In \cite[Section 2.5.d]{Katok95} it is shown that there is a semiconjugacy from a topological Markov chain to the dual dynamical system $(\T^2,A).$  From this it follows that $(\T^2,A)$ has a minimal subsystem $(X,A)$ that admits more than one invariant probability measure.  By \cite{Toms13} $C(X)\rtimes \Z$ is $\mathcal{Z}$-stable. 

\end{remark}
\begin{definition} \label{def:inner-outer} 
	Let $A$ be a C*-algebra with unique trace $\tau.$  Let $\alpha$ be an automorphism of $A$.  Then $\alpha$ preserves the unique trace and hence extends to an automorphism $\alpha$ of the von Neumann algebra $M = \pi_\tau(A)''$ where $\pi_\tau$ is the GNS representation of $A$ associated with $\tau.$ 
\begin{enumerate}
\item If $\alpha$ is an outer automorphism of $M$ then we say $\alpha$ is a \textbf{strongly outer} automorphism of $A$.   
\item If $\alpha$ is an inner automorphism of $M$ then we say $\alpha$ is a \textbf{weakly inner} automorphism of $A.$
\end{enumerate}
If $(G,\alpha,\omega)$ is a twisted action of $A$ we call the action strongly outer, weakly inner or inner if every $\alpha_t$ with $t\neq e$ has the corresponding property.
\end{definition}

We combine theorems of B\'{e}dos, and Matui and  Sato to obtain
\begin{theorem} [\textup{\cite{Bedos1993Unique, MatuiSato2014}}] \label{thm:strouter}
 Let $A$ be in $\mathcal{C}$ and $(G,\alpha,\omega)$ a strongly outer twisted action with $G$ elementary amenable.
 Then the twisted crossed product $A\rtimes_{\alpha,\omega} G$ is also in $\mathcal{C}.$
 \end{theorem}
\begin{proof}  
	The case of finite-dimensional $A$ is trivial since $A$ is a matrix algebra and has no strongly outer automorphisms, implying that $G$ is trivial and $A\rtimes_{\alpha,\omega} G \cong A \in \mathcal{C}$. 
	Hence we may assume $A$ is infinite-dimensional. 
By \cite[Theorem 1]{Bedos1993Unique}, it follows that $A\rtimes_{\alpha,\omega} G$ is simple with unique trace.\footnote{In fact his result applies to reduced crossed products of all discrete groups} 
By Corollary 4.11 and Remark 4.12 of \cite{MatuiSato2014} the crossed product $A\rtimes_{\alpha,\omega} G$ is $\mathcal{Z}$-stable and hence in $\mathcal{C}$ by \Cref{rem:equivalentCdef}.
\end{proof}
The following is clear.
\begin{lemma}\label{lem:stableC}  Let $A$ be a C*-algebra.  Then $M_n\otimes A$ is in $\mathcal{C}$ for some $n\geq1$ if and only if $M_n\otimes A$ is in $\mathcal{C}$ for all $n\geq1.$
\end{lemma}

\begin{lemma} \label{lem:innerinC} Let $A$ be in $\mathcal{C}$ and $(F,\alpha,\omega)$ a twisted inner action of a finite group on $A.$ Then $A\rtimes_{\alpha,\omega} F$ is a finite direct sum of elements of $\mathcal{C}.$
\end{lemma}
\begin{proof} By \cite[Proposition 2.3]{Eckhardt23a} we have $A\rtimes_{\alpha,\omega} F  \cong A\otimes C^*(F,\sigma)$ for some $\sigma\in Z^2(F,\T)$.  Since $F$ is finite, $C^*(F,\sigma)$ is finite dimensional. Hence the claim follows from \Cref{lem:stableC}.
\end{proof}

\begin{example} \label{example:shift-action}
	Consider the shift action $\alpha$ of a countable discrete group $G$ on the tensor product $A^{\otimes G}$, where $A \in \mathcal{C}$. Then it follows from \Cref{rem:equivalentCdef} that $A^{\otimes G} \in \mathcal{C}$ and $\pi_{\tau^{\otimes G}}(A^{\otimes G})'' \cong M^{\overline{\otimes} G}$, where $M$ is either $\mathcal{R}$ or $M_n$ for some positive integer $n$, and we write $\overline{\otimes}$ for tensor products of von Neumann algebras. 
	It follows that the induced action on the von Neumann algebra tensor product $M^{\overline{\otimes} G}$ is also the shift action. 
	
	It is well-known that the shift automorphism by an arbitrary $g \in G$ on the von Neumann algebra tensor product $M^{\overline{\otimes} G}$ is outer if and only if either $M \cong \mathcal{R}$ or $g$ has infinite order. It follows that $\alpha$ is strongly outer if either $A$ is infinite dimensional or $G$ is torsion free. In this case, we apply \Cref{thm:strouter} to conclude that $A^{\otimes G} \rtimes_\alpha G$ is also in $\mathcal{C}$, provided that $G$ is elementary amenable. 
\end{example}

\begin{remark} \label{rem:weaklyinner} Let $A\in\mathcal{C}$ and $(F,\alpha,\omega)$ a twisted action with $F$ finite. \Cref{thm:strouter} and \Cref{lem:innerinC} show if the action is either strongly outer or inner then $A\rtimes_{\alpha,\omega}F$ is a direct sum of elements of $\mathcal{C}.$  Suppose that the action is point-wise outer but not necessarily strongly outer.  In this case, the crossed product is simple.  This follows from  \cite[Theorem 1.1]{Rieffel1980Actions} (see also \cite[Theorem 10.2.2]{Barcelona18} for a less general, but shorter proof).  Moreover since $A$ is $\mathcal{Z}$-stable then so is the crossed product by a very special case of Sato's theorem \cite[Theorem 1.1]{Sato19}.  But interestingly one can not guarantee unique trace.  Indeed Chris Phillips gave examples in \cite[Section 2]{Phillips15} of outer actions of $\Z/2\Z$ on UHF algebras such that the crossed product does not have unique trace.  

In other words the presence of outer, but not strongly outer automorphisms \emph{could} lead to serious difficulties when bounding the nuclear dimension of twisted group C*-algebras.  
The next Lemma shows that we never have to worry about these types of actions for twisted group C*-algebras that appear in our analysis.
\end{remark}

\begin{lemma} \label{lem:innereq} Let $G$ be a discrete group, $\sigma\in Z^2(G,\T)$ and $H\trianglelefteq G$ a normal subgroup.  Fix $g\in G$, let $\alpha$ be the  automorphism of $H$ defined by conjugation by $g$, and let $\beta$ be the $*$-automorphism of $\C[H,\sigma]$ defined by conjugation by $u_g$. Then the following hold: 
	\begin{enumerate}
		\item \label{lem:innereq::extend} The automorphism $\beta$ extends to automorphisms of $ C^*_{\mathrm{r}}(H,\sigma) $ and $L(H,\sigma)$, which we still denote by $\beta$. 
		\item \label{lem:innereq::inner} If $L(H,\sigma)$ is a factor, then the following conditions are equivalent: 
		{
		\begin{enumerate}[ref=\alph*]
			\item \label{lem:innereq::inner::alg} $\beta$ is an inner automorphism of $\C[H,\sigma]$;
			\item \label{lem:innereq::inner::Cr} $\beta$ is an inner automorphism of $C^*_{\mathrm{r}}(H,\sigma)$;
			\item \label{lem:innereq::inner::L} $\beta$ is an inner automorphism of $L(H,\sigma)$. 
		\end{enumerate}
		}
		\item \label{lem:innereq::outer} If $L(H,\sigma)$ is a factor and the set $\{  \alpha(s)^{-1}ts: s\in H \}$ is infinite for each $t\in H$, then $\beta$ is strongly outer.
	\end{enumerate}
\end{lemma}
\begin{proof} By {\Cref{sec:twistedgroupCalgebras}}, we may suppose that $\sigma$ is normalized. Then the formula for $\beta$ is
\begin{equation*}
\beta(u_s) = u_gu_su_{g^{-1}}=\sigma(g,s)\sigma(gs,g^{-1})u_{\alpha(s)} \quad \text{ for any } s \in H \; .
\end{equation*}
Since $\beta$ preserves the canonical trace of $\C[H,\sigma]$, it extends to automorphisms of $ C^*_{\mathrm{r}}(H,\sigma) $ and $L(H,\sigma)$, which proves \eqref{lem:innereq::extend}. 

To prove \eqref{lem:innereq::inner}, we observe that clearly \eqref{lem:innereq::inner::alg} $\Rightarrow$ \eqref{lem:innereq::inner::Cr} $\Rightarrow$ \eqref{lem:innereq::inner::L}, so it suffices to show \eqref{lem:innereq::inner::L} $\Rightarrow$ \eqref{lem:innereq::inner::alg}. Suppose \eqref{lem:innereq::inner::L} holds, i.e., there is a unitary $W\in L(H,\sigma)$ such that $\beta = \text{Ad}(W).$ Decompose $W=\sum_{t \in H} x_tu_t$ with $x_t\in\C$ where convergence is in the $L^2$-norm coming from the canonical trace.  
Choose $t_0\in H$ such that $x_{t_0}\neq0$.  
For any $s\in H$, we have  
\begin{align} \label{eq:innerstuff}
	\sum_{t\in H}x_t\sigma(t,s)u_{ts}=Wu_s&=\sigma(g,s)\sigma(gs,g^{-1})u_{\alpha(s)}W\\
	&=\sigma(g,s)\sigma(gs,g^{-1})\sum_{t\in H}x_{\alpha(s)^{-1}ts}\sigma(\alpha(s),\alpha(s)^{-1}ts)u_{ts} \; , \notag
\end{align}
which implies 
\begin{equation} \label{eq:coeffcomp}
	x_{t_0}\sigma(t_0,s) = \sigma(g,s)\sigma(gs,g^{-1})x_{\alpha(s)^{-1}t_0s}\sigma(\alpha(s),\alpha(s)^{-1}t_0s)
\end{equation}
and, in particular, $|x_{t_0}|=|x_{\alpha(s)^{-1}t_0s}|$. This forces the set $X := \{ \alpha(s)^{-1}t_0s:s\in H \}$ to be finite, for otherwise $W$ would have infinite $L^2$-norm.  
Set
\begin{equation*}
	\tilde W = |X|^{-\frac{1}{2}}|x_{t_0}|^{-1}\sum_{t\in X}x_tu_t \in \C[H,\sigma]
\end{equation*}
Since $\alpha(s)^{-1}Xs=X$ for all $s\in H$, a calculation resembling \eqref{eq:innerstuff} and using \eqref{eq:coeffcomp} reveals  
\[
	\tilde{W} u_s=\sigma(g,s)\sigma(gs,g^{-1})u_{\alpha(s)} \tilde{W}  \quad \text{ for any } s \in H \; .
\]
A direct calculation (e.g., by following \cite[Lemma 4.2]{EckhardtGillaspyMcKenney2019Finite} with $A=B=L(H,\sigma)$) then shows $\tilde{W}^*\tilde{W}$ is in the center of $L(H,\sigma)$.  Since $L(H,\sigma)$ is a factor we have that $\tilde{W}^*\tilde{W} = \tau(\tilde{W}^*\tilde{W})1_{L(H,\sigma)}=1_{L(H,\sigma)}.$ Since $L(H,\sigma)$ is finite, $\tilde{W}$ is unitary and we obtain $\beta = \text{Ad}(\tilde W)$, which is the extension of an inner automorphism of $\C[H,\sigma]$.

From the proof above, we also see that if the set $= \{  \alpha(s)^{-1}ts: s\in H \}$ is infinite for each $t\in H$, then \eqref{lem:innereq::inner::L} cannot hold, that is, $\beta$ is strongly outer, which proves \eqref{lem:innereq::outer}.
\end{proof}

\begin{lemma}\label{lem:preserveinner}   Let $A$ be a finite unital C*-algebra with trivial center.  Let $n\geq 1.$   Let $\alpha$ be an automorphism of $M_n\otimes A$ that leaves $1_n\otimes A$ globally invariant. If $\alpha$ is an inner automorphism of $M_n\otimes A,$ then $\alpha|_{1_n\otimes A}$ is an inner automorphism of $A.$  
\end{lemma}
\begin{proof}  Let $\{ e_{i,j}: 1\leq i,j\leq n \}$ be matrix units for $M_n$ and $W = \sum_{i,j} e_{i,j}\otimes u_{i,j}\in M_n\otimes A$ be a unitary such that $Wx=\alpha(x)W$ for all $x\in M_n\otimes A.$  Choose indices so $u_{i,j}\neq 0.$  Since $\alpha$ leaves $A$ invariant we slightly abuse notation to obtain $u_{i,j}a=\alpha(a)u_{i,j}$ for all $a\in A$. It follows from a simple calculation (e.g., by applying \cite[Lemma 4.2]{EckhardtGillaspyMcKenney2019Finite}) that $u_{i,j}^*u_{i,j}\in A\cap A'$. Hence $u_{i,j}$ is a (non-zero) scalar multiple of a unitary in $A.$
\end{proof}
\begin{corollary}\label{lem:stablestouter} Suppose that $A$ is a unital C*-algebra and admits a unique faithful trace. Let $\alpha$ be an automorphism of $M_n\otimes A$ that leaves $1_n\otimes A$ globally invariant.  If $\alpha|_{1_n\otimes A}$ is strongly outer, then $\alpha$ is also strongly outer.
\end{corollary}
\begin{proof} Let $\tau$ be the unique trace of $A.$ Then the von Neumann algebra $\pi_\tau(A)'' = M$ is a factor.  Apply \Cref{lem:preserveinner} to $M_n\otimes M$ to obtain the conclusion.
\end{proof}

\begin{lemma}\label{lem:findexC} Let $G$ be a discrete group and $\sigma\in Z^2(G,\T).$  Suppose that $N\trianglelefteq G$ is a finite index normal subgroup such that $C^*(N,\sigma) \in \mathcal{C}.$  Then $C^*(G,\sigma)$ is a finite direct sum of elements of $\mathcal{C}.$
\end{lemma}
\begin{proof}  Set $F=G/N$ and decompose $C^*(G,\sigma)\cong C^*(N,\sigma)\rtimes_{\alpha,\omega}F$ as in \Cref{sec:groupCalgebrastwistedcp}. Let $F_1 = \{ t\in F: \alpha_t \text{ is weakly inner }  \}.$  Then $F_1$ is a normal subgroup of $F$ and by \Cref{lem:innereq} we have $\alpha_t$ is actually inner for each $t\in F_1.$ We decompose once more as in \Cref{sec:groupCalgebrastwistedcp} to obtain
\begin{equation*}
C^*(G,\sigma)\cong \left(C^*(N,\sigma) \rtimes_{\alpha,\omega} F_1\right)\rtimes_{ \alpha',\omega'} F/F_1
\end{equation*}
Notice that for each $t\in F$ we have $\alpha_t$ is inner conjugate to $\alpha'_{tF_1}|_{C^*(N,\sigma)}.$ Hence if $t\in F\smallsetminus F_1$, then $\alpha'_{tF_1}|_{C^*(N,\sigma)}$  is strongly outer.

By \cite[Proposition 2.3]{Eckhardt23a}, we have an isomorphism $\pi:  C^*(N,\sigma) \rtimes_{\alpha,\omega} F_1\to C^*(N,\sigma) \otimes C$ for some finite dimensional C*-algebra $C$ where $\pi(a) = a\otimes 1_C$ for all $a\in C^*(N,\sigma).$ We then have
\begin{equation*}
C^*(G,\sigma)\cong \left( C^*(N,\sigma) \otimes C\right)\rtimes_{\alpha',\omega'} F/F_1
\end{equation*}
where each automorphism $\alpha_{tF_1}'$ leaves $C^*(N,\sigma) \cong C^*(N,\sigma) \otimes 1_C$ globally invariant.
 
Let $P = \{ p_1,...,p_n \}$ be the minimal central projections of $C^*(N,\sigma) \otimes C.$  Then $\alpha'$ restricts to a true action of $F/F_1$ on $P.$ By considering each orbit separately and taking the direct sum we may, without loss of generality, assume that $F/F_1$ acts transitively on $P.$   Let $H\leq F/F_1$ be the stabilizer of $p_1.$ We have $p_1C^*(N,\sigma) \otimes C\cong C^*(N,\sigma)\otimes M_k$ for some $k.$   By \cite[Proposition 2.2]{Eckhardt23a} we have
\begin{equation*}
\left( C^*(N,\sigma) \otimes C\right)\rtimes_{\alpha',\omega'} F/F_1\cong M_n\otimes (C^*(N,\sigma)\otimes M_k \rtimes_{\alpha',\omega'}H)
\end{equation*}

By \Cref{lem:stablestouter} for each non-trivial $t\in H$ the automorphism $\alpha_t'$ restricted to $C^*(N,\sigma)\otimes M_k$ is strongly outer.  Then \Cref{thm:strouter} and \Cref{lem:stableC} show that $M_n\otimes (C^*(N,\sigma)\otimes M_k \rtimes_{\alpha',\omega'}H)
\in \mathcal{C}.$
\end{proof}

\section{Nuclear dimension bounds} \label{sec:FND}
In the three subsections below, we prove our main theorems on finite nuclear dimension for (twisted) group C*-algebras. 
As indicated in the introduction, 
these theorems rely on a recent result of Hirshberg and the second named author in \cite{HirshbergWuActions}, which is about bounding the nuclear dimension of crossed products of $G$-$X$-C*-algebras. In \Cref{thm:HW} below, we give a generalization of said result to the case of twisted crossed products and show how this generalization follows easily from the Packer-Raeburn stabilization trick described in \Cref{sec:stabilization}. 
The statement of \Cref{thm:HW} makes use of the following concepts that we treat as blackboxes, both being dimensions associated to an action $\beta \colon G \curvearrowright X$ of a discrete group $G$ on a locally compact Hausdorff space $X$: 
\begin{enumerate}
	\item the long thin covering dimension $\dimltc(\beta)$,  
	introduced in \cite[Section~5]{HirshbergWuActions}, and 
	\item the asymptotic dimension of the coarse orbit space of $\beta$, denoted by $\mathrm{asdim}(X, \mathcal{E}_\beta)$ and introduced in \cite[Section~3]{HirshbergWuActions}. 
\end{enumerate}
For the purpose of proving our main theorem, it suffices to know that in the setting we care about, we have the following explicit upper bounds for these dimensions: if $G$ is a finitely generated virtually nilpotent group and $k$ is the degree of its polynomial growth function, then it follows from \cite[Theorems~7.4, 7.5, and~8.1, as well as Remark~A.3 and Lemma~A.4]{HirshbergWuActions} that we have 
\begin{equation} \label{eq:bounds}
	\mathrm{asdim}(X, \mathcal{E}_\beta) + 1 \leq 3^k
	\quad \text{ and } \quad 
	\dimltc(\beta)  + 1 \leq  3^k (\dim(X) + 1 ) \; .
\end{equation}

\begin{theorem} \label{thm:HW}
	Let $(G,A,\alpha,\omega, X, \beta, \iota)$ be a twisted $G$-$X$-C*-algebra, where $G$ is discrete and $X$ is a locally compact Hausdorff space. 
	Then we have 
	\[
		\dimnuc(A \rtimes_{\alpha, \omega, \mathrm{max}} G) +1 \leq (\mathrm{asdim}(X, \mathcal{E}_\beta) +1) \cdot
		(\dimltc(\beta) +1)  
		\cdot (d_{\mathrm{stab}} + 1) 
		\; ,
	\]
	where, using notations from \Cref{sec:G-X-algebras}, we define 
	\[
		d_{\mathrm{stab}} = \sup \left\{ \dimnuc \left(A_x \rtimes_{\alpha_{|H}^x, \omega_{|H}^x, \mathrm{max} } H \right) \colon x \in X, H 
		\leq G_x \right\} 
		\; .
	\]
\end{theorem}

\begin{proof}
	It follows from the construction in \Cref{sec:stabilization} that $(G, A  \otimes \mK, \alpha, \omega \otimes 1)$ is exterior equivalent to a genuine action $\alpha^\omega \colon G \curvearrowright A  \otimes \mK$. 
	It follows from \Cref{sec:G-X-algebras} that for any $x \in X$ that is fixed by a subgroup $H$ of $G$, the induced twisted dynamical system $(H, A_x \otimes \mK, \alpha^x_{|H}, \omega^x_{|H} \otimes 1)$ is exterior equivalent to the genuine action $(\alpha^\omega)^x_{|H} \colon H \curvearrowright A_x \otimes \mK$. 
	Then it follows from \Cref{sec:twisted-croe} that we have the following $*$-isomorphisms: 
	\begin{multline*}
		(A_x \rtimes_{\alpha^x_{|H}, \omega^x_{|H} , \mathrm{max}} H ) \otimes \mK
		\cong (A_x \otimes \mK) \rtimes_{\alpha^x_{|H}, \omega^x_{|H} \otimes 1, \mathrm{max}} H 	\cong (A_x \otimes \mK) \rtimes_{(\alpha^\omega)^x_{|H}, \mathrm{max}} H  
	\end{multline*}
	and 
	\[
		(A \rtimes_{\alpha, \omega , \mathrm{max}} G ) \otimes \mK
		\cong (A \otimes \mK) \rtimes_{\alpha, \omega \otimes 1, \mathrm{max}} G 	\cong (A \otimes \mK) \rtimes_{\alpha^\omega, \mathrm{max}} G  \; .
	\]
	Since by \Cref{lemma:dimnuc-basic}\eqref{lemma:dimnuc-basic::stabilization}, 
	taking stabilization preserves the nuclear dimension, 
	the desired inequality then follows from applying \cite[Theorem~10.2]{HirshbergWuActions} (which is precisely \Cref{thm:HW} for untwisted crossed products) to the genuine $G$-$X$-C*-algebra $(G,A \otimes \mK,\alpha^\omega, X, \beta, \iota \otimes 1)$. 
\end{proof}

\subsection{Virtually polycyclic groups}
We will show in this subsection that the group C*-algebra of a virtually polycyclic group has finite nuclear dimension. 
We start with a more refined discussion of the construction in \Cref{ex:beyond-central-extension}. We briefly recall the latter here (with the notation $H$ replacing $G$): Given discrete groups $H$, $N$, and $M$ with $N$ and $M$ normal in $H$ and $M \leq Z(N)$, we can view $C^*(N)$ as a twisted $H/N$-$\widehat{M}$-C*-algebra such that $C^*(H)$ is identified with the associated twisted crossed product $C^*(N) \rtimes_{\alpha,\omega, \mathrm{max}} H/N $, where we denote the twisted action of $H/N$ on $C^*(N)$ by $(\alpha,\omega)$. 

\begin{lemma} \label{lemma:beyond-central-extension-fixed-char}
	Let $H$, $N$, $M$ and $(\alpha,\omega)$ be as above. 
	Let $\chi$ be a character on $M$ that is $H$-invariant in the sense that $\chi \circ \Ad_h = \chi$ for any $h \in H$. 
	This induces a twisted action $(\alpha_\chi, \omega_\chi)$ of $H/N$ on $C^*(N)_\chi$. 
	Let $L = \ker (\chi) \leq M$ and keep the notation $\chi$ for the induced character on $M / L$. Then the following hold: 
	\begin{enumerate}
		\item \label{lemma:beyond-central-extension-fixed-char::normal} The kernel $L$ is a normal subgroup in $H$. 
		\item \label{lemma:beyond-central-extension-fixed-char::central} The quotient $M / L$ is central in $H/L$. 
		\item \label{lemma:beyond-central-extension-fixed-char::quotient} The twisted crossed product $C^*(N)_\chi \rtimes_{\alpha_\chi,\omega_\chi,\mathrm{max}} H/N$ is a quotient of $C^*(H/L)$. 
		\item \label{lemma:beyond-central-extension-fixed-char::twisted} The twisted crossed product $C^*(N)_\chi \rtimes_{\alpha_\chi,\omega_\chi,\mathrm{max}} H/N$ is isomorphic to 
		$C^*(H/M, \sigma_\chi)$, where $\sigma_\chi$ is associated to the induced character $\chi$ on $M/L$ as in \Cref{ex:central-extension}.
	\end{enumerate}
\end{lemma}

\begin{proof}
	Statement~\eqref{lemma:beyond-central-extension-fixed-char::normal} follows directly from the $H$-invariance of $\chi$. 
	Statement~\eqref{lemma:beyond-central-extension-fixed-char::central} follows from a direct computation: 	
	for any $h \in H$ and any $m \in M$, we have
	\begin{equation*}
	\chi([h, m]) = \chi(h^{-1}m^{-1}h m)=\chi(h^{-1}m^{-1}h)\chi(m)=\chi(m^{-1})\chi(m) = 1 \; ,
	\end{equation*}
	that is, $[h,m] \in L$. 
	
	To prove \eqref{lemma:beyond-central-extension-fixed-char::quotient} and \eqref{lemma:beyond-central-extension-fixed-char::twisted}, we observe there is a canonical $*$-isomorphism $C^*(N)_\chi \cong C^*(N/L)_\chi$.  
	We recall the choice of an identity preserving lift $H/N \to H$ in \Cref{ex:beyond-central-extension} and use it to construct an identity preserving lift $c'$ from the composition $(H/L)/(N/L) \xrightarrow{\cong} H/N \xrightarrow{c} H \twoheadrightarrow H / L$, which leads to a $*$-isomorphism 
	\begin{equation} \label{lemma:beyond-central-extension-fixed-char::eq:mod-L}
	C^*(N)_\chi\rtimes_{\alpha_\chi , \omega_\chi,\mathrm{max}} H/N \cong C^*(N/L)_\chi\rtimes_{\alpha'_\chi , \omega'_\chi,\mathrm{max} } (H/L)/(N/L) 
	\end{equation}
	where 
	\begin{equation*} 
	\omega'_\chi(s,t)=\left( u_{c'(s)}u_{c'(t)}u_{c'(st)}^{*} \right)_\chi \; \textup{ and } \; (\alpha'_\chi)_s=\text{Ad}\left(\left(u_{c'(s)}\right)_\chi \right) 
	\quad \text{ for any } s, t \in (H/L)/(N/L) \; .
	\end{equation*}
	The same formulas without the restrictions to $\chi$ give rise to a twisted action $(\alpha' , \omega')$ of $(H/L)/(N/L)$ on $C^*(N/L)$.  
	Comparing this with \eqref{eq:tgastcp-action}, we obtain a $*$-isomorphism 
	\begin{equation} \label{lemma:beyond-central-extension-fixed-char::eq:H-N}
	C^*(N/L) \rtimes_{\alpha',\omega',\mathrm{max}} (H/L)/(N/L) \cong C^*(H/L) \; .
	\end{equation}
	Moreover, thanks to \eqref{lemma:beyond-central-extension-fixed-char::central}, we have a central extension $1 \to M / L \to H / L \to H / M \to 1$. 
	Choosing some identity preserving lift $d \colon H/M \to H/L$ and applying the construction in \Cref{ex:central-extension} to this central extension, we can view $C^*(H / L)$ as an $\widehat{M/L}$-C*-algebra, such that 
	\begin{equation} \label{lemma:beyond-central-extension-fixed-char::eq:chi-M-twisted}
	C^*(H / L)_\chi \cong C^*(H/ M,\sigma_\chi) \quad \text{ with } \sigma_\chi ([s], [t]) = \chi \left(d([s])d([t])d([st])^{-1} \right)  \text{ for any } s,t \in H  \; . 
	\end{equation}
	Note that the $*$-isomorphism in \eqref{lemma:beyond-central-extension-fixed-char::eq:H-N} respects the $\widehat{M/L}$-C*-algebra structures, so it induces a $*$-isomorphism
	\begin{equation} \label{lemma:beyond-central-extension-fixed-char::eq:H-N-chi}
	C^*(N/L)_\chi\rtimes_{\alpha'_\chi , \omega'_\chi ,\mathrm{max}} (H/L)/(N/L) \cong C^*(H / L)_\chi  \; .
	\end{equation}
	Combining \eqref{lemma:beyond-central-extension-fixed-char::eq:mod-L} with \eqref{lemma:beyond-central-extension-fixed-char::eq:H-N-chi} proves \eqref{lemma:beyond-central-extension-fixed-char::quotient}, while a further combination with \eqref{lemma:beyond-central-extension-fixed-char::eq:chi-M-twisted} proves \eqref{lemma:beyond-central-extension-fixed-char::twisted}. 
\end{proof}

The following lemma relates our discussion of strongly outer actions in \Cref{sec:mathcalC} to the setting of polycyclic groups, where the Fitting subgroup (see \Cref{sec:Polycyclicgroups}) plays a vital role. 

\begin{lemma}\label{lem:strongouterfitting} Let $K$ be a polycyclic group such that its Fitting subgroup $N=\fitt(K)$ is torsion free and that $K/N$ is abelian. Let $\chi\in \widehat{Z(N)}$ be a faithful character and   suppose that $k\in K\smallsetminus N$ such that $\chi(x)=\chi(kxk^{-1})$ for all $x\in Z(N).$  
Set $H = \la  N,k \ra.$  Then 
$Z(N) \leq Z(H)$. 
Moreover, let $\sigma_\chi\in Z^2(H,\T)$ be defined as in \Cref{ex:central-extension} and 
let $\beta$ be the automorphism of $C^*(N,\sigma_\chi)$ given by conjugation by $u_k \in C^*(H,\sigma_\chi)$. 
Then $\beta$ is strongly outer.
\end{lemma}
\begin{proof}  
	Since $\chi$ is faithful and $\chi([k, x]) = \chi(k^{-1}xk)^{-1} \chi(x) = 1$ for all $x\in Z(N)$, we have $Z(N) \leq Z(H)$. 
	Now for any $x\in K$, let $\overline{x}$ denote the image of $x$ in $K/Z(N).$ 
By \Cref{lem:innereq}, it suffices  to show that the set $\{ \overline{k}s^{-1}\overline{k}^{-1}as:s\in N/Z(N) \}$ is infinite for all $a\in N/Z(N)$.  We prove the equivalent statement that 
$\{ sa\overline{k}s^{-1}:s\in N/Z(N) \}$ is infinite for all $a\in N/Z(N)$. 

Suppose the contrary for some $a\in N/Z(N)$. 
Then observe that $[a\overline{k}, s] = 1$ for any $s\in N/Z(N)$, since it follows from the finiteness of $\{ s^n a\overline{k}s^{-n}: n \in \mathbb{Z} \}$ that there are integers $n<m$ such that $s^na\overline{k}s^{-n}=s^ma\overline{k}s^{-m}$, which can be rearranged to $[a\overline{k}, s^{m-n}] = 1$ and thus implies $[a\overline{k}, s] = 1$ by \Cref{lem:uniqueroot}. 
Hence 
for any $b\in N$ with $a=\overline{b}$, we have
\begin{equation*} \label{eq:commincenter}
[bk,x]\in Z(N) \textup{ for all }x\in N.
\end{equation*}
It follows that $H/Z(H)$ is a quotient of $N/Z(N)\times \Z$ via mapping a generator of $\mathbb{Z}$ to the class of $bk$. 
Hence $H/Z(H)$ and thus also $H$ are nilpotent, too. 
 
 Since $K/N$ is abelian, $H$ is also a normal subgroup of $K$ so $bk\in H\leq \fitt(K)=N.$ Since $b\in N$, we have $k\in N$, a contradiction. 
\end{proof}

\begin{theorem}\label{thm:finitend}  There is an increasing function $f:\Z^{\geq 0}\rightarrow \Z^{\geq0}$ such that for every virtually polycyclic group $G$ we have $\dimnuc(C^*(G)) \leq f(\hirsch(G)) < \infty$.  
\end{theorem}
\begin{proof}  We proceed by induction on the Hirsch length. The starting case of $\hirsch(G)=0$ is equivalent to assuming $G$ is finite, in which case $C^*(G)$ has nuclear dimension 0 so we set $f(0)=0$.  
	The next case of $\hirsch(G)=1$ is equivalent to assuming $G$ is virtually cyclic, that is, it contains a finite index normal subgroup isomorphic to $\mathbb{Z}$. In this case we may decompose $C^*(G)$ as a twisted crossed product $C^*(\mathbb{Z}) \rtimes_{\alpha, \omega} F$ as in \Cref{sec:groupCalgebrastwistedcp} where $F := G / \mathbb{Z}$ is finite. Since $C^*(\mathbb{Z}) \cong C(\mathbb{T})$ by the Gelfand duality and the automorphism group of $\mathbb{Z}$ is $\mathbb{Z}/2\mathbb{Z}$, with its induced action on $\mathbb{T}$ given by the flip about the real axis, whose quotient space can be canonically identified with $[0,\pi]$, we see that $C^*(G)$ is a twisted $[0,\pi]$-C*-algebra with finite-dimensional fibers. It follows from  \Cref{lemma:dimnuc-basic}\eqref{lemma:dimnuc-basic::field} that $\dimnuc(C^*(G)) \leq 1$. So we may set $f(1)=1$.  
	
Now for any integer $n > 1$, 
supposing 
$f(n-1)$ has been defined and satisfies the conclusion of the theorem, we define  
\begin{equation} \label{thm:finitend::eq:recursive}
	f(n)=9^{n} (n + 1) (f(n-1) + 1) - 1  \; .
\end{equation}
Clearly we have $f(n) > f(n-1)$. It remains to show that for any virtually polycyclic group $G$ with $\hirsch(G) = n$, we have $\dimnuc(C^*(G)) \leq f(n).$

To this end, as mentioned in \Cref{sec:Polycyclicgroups}, 
there are characteristic subgroups $N\trianglelefteq K\trianglelefteq G $ such that
\begin{enumerate}
\item $G/K$ is finite, 
\item $N$ is the Fitting subgroup of $K$, and
\item $K$ is torsion free and $K/N$ is free abelian. 
\end{enumerate}

Following \Cref{ex:beyond-central-extension}, we can identify $C^*(N)$ as a twisted $G/N$-$\widehat{Z(N)}$-C*-algebra in the sense of \Cref{sec:G-X-algebras} so that 
\begin{equation*}
C^*(G) \cong C^*(N) \rtimes_{\alpha,\omega} G/N \; ,
\end{equation*}
to which we can apply \Cref{thm:HW},  
after some preparation:  
Since $Z(N)$ is finitely generated abelian, it follows from \eqref{eq:Hirsch-Pontryagin} that 
$\dim(\widehat{Z(N)}) = \hirsch(Z(N))\leq \hirsch(G)$.    
Since $G/N$ is virtually abelian the degree of its polynomial growth function is equal to $\hirsch(G/N)$ and hence less than or equal to $\hirsch(G).$  At this point we combine the  
 inequalities \eqref{eq:bounds} and \Cref{thm:HW} to see that
\begin{equation*}
\dimnuc(C^*(G)) + 1 \leq 3^{\hirsch(G)}\cdot 3^{\hirsch(G)}(\hirsch(G) + 1)\cdot (d_{\mathrm{stab}} + 1) \; ,
\end{equation*}
where, keeping the notation $(\alpha,\omega)$ when we restrict it to a subgroup $H/N$ of $G/N$ and a fiber $C^*(N)_\chi$ fixed by $H/N$, we write 
\[
	d_{\mathrm{stab}} := \sup \left\{ \dimnuc \left(C^*(N)_\chi\rtimes_{\alpha,\omega} H/N \right) \colon \chi \in \widehat{Z(N)}, N\leq H\leq G_\chi \right\} 
\; ,
\]
where $G_\chi := \{ h\in G: \chi\circ \text{Ad}(h) = \chi  \},$  is the stabilizer of $\chi$ under the conjugation action.
In view of \eqref{thm:finitend::eq:recursive}, to complete the proof, it suffices to show that $d_{\mathrm{stab}} \leq f(\hirsch(G) - 1)$, 
that is, for any fixed $\chi\in\widehat{Z(N)}$ and for any fixed subgroup $H \leq G$ satisfying $N \leq H \leq G_\chi$, we have 
\[
	\dimnuc \left(C^*(N)_\chi\rtimes_{\alpha,\omega} H/N \right) \leq f(\hirsch(G) - 1) \; .
\]

To this end, 
we discuss two cases:

\begin{enumerate}
	\item Suppose first that $\text{ker}(\chi)$ is not trivial. By \Cref{lemma:beyond-central-extension-fixed-char}\eqref{lemma:beyond-central-extension-fixed-char::quotient}, we may view $C^*(N)_\chi\rtimes_{\alpha,\omega} H/N$ as a quotient of the group C*-algebra $C^*(H/\text{ker}(\chi)).$ Since $N$ is torsion free we have $\hirsch(\ker(\chi))\geq 1$, hence $\hirsch(H/\text{ker}(\chi)) = \hirsch(H) - \hirsch(\text{ker}(\chi))\leq \hirsch(G)-1.$ Hence
	\begin{equation*}
	\dimnuc(C^*(N)_\chi\rtimes_{\alpha,\omega} H/N)\leq \dimnuc(C^*(H/\text{ker}(\chi))) \leq f(\hirsch(G)-1).
	\end{equation*}
	
	\item Next suppose that $\text{ker}(\chi)$ is trivial. 
	By \Cref{lemma:beyond-central-extension-fixed-char}\eqref{lemma:beyond-central-extension-fixed-char::twisted}, 
	we have
	\begin{equation*}
	C^*(N)_\chi\rtimes_{\alpha,\omega} H/N \cong C^*(H/Z(N),\sigma_\chi) \; .
	\end{equation*}
	Let $N\leq H'\leq H$ so that $H'/N = K/N \cap H/N.$  By \Cref{lem:strongouterfitting}, for each $t\in H'/N$ the automorphism $\alpha_t$ is strongly outer. By \Cref{thm:nilpotenttwist} the fiber $C^*(N)_\chi\in\mathcal{C}.$ Since $H'/N$ is abelian by \Cref{thm:strouter} we have
	\begin{equation*}
	C^*(N)_\chi \rtimes_{\alpha,\omega} H'/N \cong C^*(H'/Z(N),\sigma_\chi) \in\mathcal{C}.
	\end{equation*}
	Since $H'$ is finite index in $H$, by \Cref{lem:findexC},  $C^*(H/Z(N),\sigma_\chi)$ is isomorphic to a direct sum of elements of $\mathcal{C}$, whence
	\[
		\dimnuc(C^*(N)_\chi\rtimes_{\alpha,\omega} H/N) = \dimnuc(C^*(H/Z(N),\sigma_\chi))\leq 1 = f(1) \leq f(\hirsch(G) - 1) \; .
	\]
\end{enumerate}
Combining the two cases, we obtain the desired inequality. 
\end{proof}

\begin{remark} 
	If one keeps track of the recursive formula \eqref{thm:finitend::eq:recursive}, one can derive an explicit formula for the bounding function $f$ in \Cref{thm:finitend}: 
	\[
	f(n) = 
	\begin{cases}
	0 , & n = 0 \\
	(n + 1)! \cdot  9^{\frac{(n+2)(n-1)}{2}}  - 1 , & n = 1,2,3, \ldots 
	\end{cases}
	\; .
	\]
	However, beyond the first few values, this bound appears to be so large that its precise value is of no practical use. This is why we do not emphasize this exact formula in the statement of \Cref{thm:finitend}. 
\end{remark}

\begin{corollary}\label{cor:indlimfdr} Let $G$ be an inductive limit of virtually polycyclic groups $G_n$ such that $\varliminf_n \hirsch(G_n)<\infty$.  Then $\dimnuc(C^*(G))<\infty$
\end{corollary}

\begin{proof}
	It follows from \Cref{lemma:dimnuc-basic}\eqref{lemma:dimnuc-basic::limit} and the integrality of the Hirsch length that 
	\[
		\displaystyle \dimnuc(C^*(G)) = \dimnuc(\varinjlim_n C^*(G_n)) \leq \varliminf_n \dimnuc(C^*(G_n)) \leq f(\varliminf_n \hirsch(G_n)) < \infty \; ,
	\]
	where $f$ is taken from \Cref{thm:finitend}. 
\end{proof}

We say a C*-algebra is \emph{classifiable} if it is simple, nuclear, $\mathcal{Z}$-stable and satisfies the universal coefficient theorem (UCT) of Rosenberg and Schochet  \cite{RosenbergSchochet1987UCT}.  By  \cite{Gong20,Carrion23}  any two classifiable C*-algebras with the same Elliott invariant are isomorphic.
\begin{corollary} Let $G$ be a virtually polycyclic group and $J$ a maximal ideal of $C^*(G).$ If $A=C^*(G)/J$ satisfies the UCT, then $A$ is classifiable.
\end{corollary}
\begin{proof}
	The simplicity of $A$ follows from the maximality of $J$, and the separability of $A$ follows from the countability of $G$. By \Cref{thm:finitend} and \Cref{lemma:dimnuc-basic}\eqref{lemma:dimnuc-basic::extension}, $A$ has finite nuclear dimension, which implies it is nuclear and $\mathcal{Z}$-stable by \cite{Winter2012nuclear}. 
\end{proof}

\begin{corollary}  Let $G$ be a virtually polycyclic group and $N\trianglelefteq G$ a virtually abelian subgroup. If $\tau$ is a trace on $G$ such that $\tau(G\smallsetminus N) = \{0\}$ and $C^*_\tau(G)$\textemdash the C*-algebra generated by the GNS representation of $\tau$\textemdash is simple, then $C^*_\tau(G)$ is classifiable.
 \end{corollary}
 \begin{proof} The C*-algebra $C^*_\tau(G)$ has finite nuclear dimension by \Cref{thm:finitend} and satisfies the UCT by \cite[Theorem 2.6]{Eckhardt23a}.
 \end{proof}

Using \Cref{lem:centralextension}, we can extend \Cref{thm:finitend} to cover all twisted group C*-algebras of virtually polycyclic groups.

\begin{theorem} \label{thm:twistednucdim}  Let $G$ be a virtually polycyclic group and $\sigma\in Z^2(G,\T).$ Then  
\begin{equation*}
\dimnuc(C^*(G,\sigma))\leq f(\hirsch(H_2(G))+\hirsch(G))<\infty \; ,
\end{equation*}
where $f$ is the function from {\Cref{thm:finitend}}. 
\end{theorem}
\begin{proof} Since $G$ is virtually polycyclic, the second homology group $H_2(G)$ is finitely generated by  \cite[Theorem 5.4]{Baumslag1981Computable}.  Let $E$ be the group constructed in \Cref{lem:centralextension}.  Since $G$ and $H_2(G)$ are both virtually polycyclic, the group $E$ is also virtually polycyclic with $\hirsch(E)=\hirsch(H_2(G))+\hirsch(G).$
By \Cref{thm:finitend} we have $\dimnuc(C^*(E))\leq f(\hirsch(E)).$  By \Cref{lem:centralextension}, $C^*(G,\omega)$ is a quotient of $C^*(E)$, hence $\dimnuc(C^*(G,\omega))\leq \dimnuc(C^*(E))$ by  \Cref{lemma:dimnuc-basic}\eqref{lemma:dimnuc-basic::extension}. 
\end{proof}

\begin{corollary} \label{cor:twistedclassifiable} Let $G$ be an infinite virtually polycyclic group and $\sigma\in Z^2(G,\sigma).$  If $C^*(G,\sigma)$ is simple then it is classifiable.
\end{corollary}
\begin{proof} By \cite[Theorem 1.1]{Barlak17},  $C^*(G,\sigma)$ satisfies the UCT and by 
\Cref{thm:twistednucdim} it has finite nuclear dimension.
\end{proof}
 \begin{remark} For \emph{nilpotent} groups, Kleppner's condition provides a straightforward test (see \cite[Theorem 1.7]{Packer89twisted}) for simplicity of $C^*(G,\sigma).$  In general it can be difficult to decide simplicity of a twisted group C*-algebra.  The reason that the nilpotent case is so much easier boils down to the fact that their group C*-algebras have a $T_1$ primitive ideal space.
 \end{remark}  
 \subsection{A class of non-residually finite examples} In \cite[Page 349]{Hall61} P. Hall constructs a finitely generated group $G$ that is nilpotent by cyclic (and thus solvable) but is not Hopfian. This implies that $G$ is not residually finite and, in particular, not virtually polycyclic.  
Hall's group is a starting point for many counterexamples in group theory and can be a good test case for extending results beyond the residually finite solvable world.

As far as we know all finitely generated examples of groups previously known to yield finite nuclear dimension are residually finite.  Our modest goal in this section is to show that with $G$ as defined above, $C^*(G)$ has finite nuclear dimension to (i) point out that residual finiteness is not necessary for finite nuclear dimension of finitely generated groups and (ii) to highlight that our methods can be used outside of the virtually polycyclic case. We also observe that $C^*(G)$ is not strongly quasidiagonal and therefore has infinite decomposition rank, echoing \Cref{conj:dr}.

We recall Hall's construction. Let $p$ be a prime number and consider the Heisenberg group over the ring $\Z[\tfrac{1}{p}]$
\begin{equation*}
\mathbb{H}_3(\Z[\tfrac{1}{p}]) = \left\{  \left[ \begin{array}{ccc} 1 & a & c\\ 0 & 1 & b\\ 0 & 0 & 1  \end{array}  \right] : a,b,c \in \Z[\tfrac{1}{p}] \right\}
\end{equation*}
Then $\mathbb{H}_3(\Z[\tfrac{1}{p}])$ is  nilpotent with center isomorphic to $\Z[\tfrac{1}{p}].$   Let $\beta$ be the automorphism 
\begin{equation*}
\beta\left(  \left[ \begin{array}{ccc} 1 & a & c\\ 0 & 1 & b\\ 0 & 0 & 1  \end{array}  \right]  \right) = \left[ \begin{array}{ccc} 1 & pa & c\\ 0 & 1 & \frac{1}{p}b\\ 0 & 0 & 1  \end{array}  \right]
\end{equation*}
Set $H = \mathbb{H}_3(\Z[\tfrac{1}{p}])\rtimes_\beta\mathbb{Z}.$  Let $u\in H$ be the element implementing $\beta.$  Then $H$ is generated by $u$ and the following two elements
\begin{equation*}
\left[  \begin{array}{ccc} 1 & 1 & 0\\ 0 & 1 & 0\\ 0 & 0 & 1  \end{array} \right],\left[  \begin{array}{ccc} 1 & 0 & 0\\ 0 & 1 & 1\\ 0 & 0 & 1  \end{array} \right].
\end{equation*} 
Let $z = \left[  \begin{array}{ccc} 1 & 0 & 1\\ 0 & 1 & 0\\ 0 & 0 & 1  \end{array} \right].$ P. Hall shows \cite[Page 349]{Hall61} that $G=H/\la  z \ra$ is not Hopfian and in particular not residually finite.  

\begin{theorem}\label{thm:nonRFexample}
Let $H$ and $G$ be the groups defined above.  Then 
\begin{equation*}
\dimnuc(C^*(H)), \dimnuc(C^*(G))\leq 53. 
\end{equation*}
Moreover neither $C^*(H)$ nor $C^*(G)$ are strongly quasidiagonal hence they have infinite decomposition rank. Both $H$ and $G$ have finite Hirsch length as defined by Hillman.
\end{theorem}
\begin{proof}  Again, following \cite{PackerRaeburn1992structure} we decompose $C^*(H)$ as a continuous field over $\widehat{\Z[\tfrac{1}{p}]}$, where the fiber over $\chi$ is the twisted group C*-algebra $C^*(H/Z(H),\sigma_\chi).$  We further decompose $C^*(H/Z(H),\sigma_\chi)$ as a crossed product
\begin{equation*}
C^*(H/Z(H),\sigma_\chi)\cong  C^*(\Z[\tfrac{1}{p}]^2,\sigma_\chi)\rtimes_\beta \Z
\end{equation*}
where we abuse notation and have $\beta(u_{(x,y)}) = u_{(px,\frac{y}{p})}.$ The analysis splits into two cases.
\\\\   
\textbf{Case 1:}  Suppose the kernel of $\chi$ is not cyclic. By Lemmas 5.6 and 5.22 of \cite{Eckhardt23} the C*-algebra $C^*(\Z[\tfrac{1}{p}]^2,\sigma_\chi)$ is a $C(X)$-algebra with $X$ homeomorphic to the 2 dimensional space $\widehat{\Z[\tfrac{1}{p}]}^2$, finite dimensional fibers and $\beta$ restricts to an action on $C(X).$ Then by  \cite[Theorem 10.2]{HirshbergWuActions} we have 
\begin{equation*}
\dimnuc(C^*(\Z[\tfrac{1}{p}]^2,\sigma_\chi)\rtimes_\beta \Z) \leq 3^1\cdot3^1(2 + 1)\cdot(0+1) - 1 = 26.
\end{equation*}
\textbf{Case 2:} Suppose the kernel of $\chi$ is cyclic. Then by the proof of \cite[Theorem 5.23]{Eckhardt23} $C^*(\Z[\tfrac{1}{p}]^2,\sigma_\chi)$ is simple with unique trace. 
By \cite[Theorem 2.2]{Carrion11}, $\dimnuc(C^*(\mathbb{H}_3(\Z)))\leq 20.$  Since $\mathbb{H}_3(\Z[\tfrac{1}{p}])$ is a direct limit of integral Heisenberg groups we have $\dimnuc(\mathbb{H}_3(\Z[\tfrac{1}{p}]))\leq 20.$ Since $C^*(\Z[\tfrac{1}{p}]^2,\sigma_\chi)$ is a quotient of $C^*(\mathbb{H}_3(\Z[\tfrac{1}{p}]))$ we have $\dimnuc(C^*(\Z[\tfrac{1}{p}]^2,\sigma_\chi))\leq 20.$ 
By \cite[Lemma 5.20]{Eckhardt23}, $\beta$ is a strongly outer automorphism and therefore by \Cref{thm:strouter} we have
\begin{equation*}
\dimnuc(C^*(\Z[\tfrac{1}{p}]^2,\sigma_\chi)\rtimes_\beta \Z)\leq 1.
\end{equation*}
For the nuclear dimension estimate, we use 
\Cref{lemma:dimnuc-basic}\eqref{lemma:dimnuc-basic::field} 
to bound the nuclear dimension of $C^*(H)$ in terms of its fibers.  
Since $\widehat{\Z[\tfrac{1}{p}]}$ is 1 dimensional we have
\begin{equation*}
\dimnuc(C^*(H)) \leq (1+1)(26+1) - 1 = 53.
\end{equation*}
Next we show infinite decomposition rank.
Consider the group $\Z[\tfrac{1}{p}]\rtimes_\alpha \Z$  where $\alpha(x)=px.$  We claim that $C^*(\Z[\tfrac{1}{p}]\rtimes_\alpha \Z)$ is not strongly quasidiagonal.  Since the groups $H$ and $G$ both quotient onto $\Z[\tfrac{1}{p}]\rtimes_\alpha \Z$ this will complete the proof.  Note that 
\begin{equation*}
\widehat{\Z[\tfrac{1}{p}]} \cong \{ (a_n)_{n=0}^\infty\in\T^\infty : a_{n+1}^p=a_n  \}
\end{equation*}
and under this identification the dual automorphism $\hat\alpha$ is given by
\begin{equation} \label{eq:alphadef}
\hat\alpha(a_0,a_1,a_2...) =(a_0^p, a_0, a_1, ...), \quad \hat\alpha^{-1}(a_0,a_1,a_2...) =(a_1,a_2, ...)
\end{equation}
For any $\lambda\in\T$ there is a $p$th root of $\lambda$ whose imaginary part is less than or equal to 0.
Therefore there is an element $\chi = (a_n)_{n=1}^\infty\in\widehat{\Z[\tfrac{1}{p}]}$ such that $a_0=1$ and for each $n\geq1$ we have $\text{Re}(a_n)\leq0.$ Then by \eqref{eq:alphadef} we see that
\begin{equation*}
\overline{\{\hat\alpha^n(\chi): n<0\}}\cap\overline{\{\hat\alpha^n(\chi): n\geq0\}} = \emptyset.
\end{equation*}
Hence the representation of $C^*(\Z[\tfrac{1}{p}]\rtimes_\alpha \Z)$ induced from the closed orbit of $\chi$ generates a C*-algebra with a proper isometry by \cite[Theorem 9]{Pimsner83} hence the C*-algebra is not strongly quasidiagonal.

Finally notice that $\Z[\tfrac{1}{p}]$ has Hirsch length 1 since $\Z[\tfrac{1}{p}]/\Z$ is locally finite.  Since 
\begin{equation*}
H \cong (\Z[\tfrac{1}{p}]^2\rtimes \Z[\tfrac{1}{p}])\rtimes \Z 
\end{equation*}
it follows that $\hirsch(H) = 4$ and $\hirsch(G)=3.$
\end{proof}

\subsection{Wreath products of finite groups by virtually nilpotent groups} \label{sec:wr-finite}
In this subsection, we verify finite nuclear dimension for yet another class of finitely generated elementary amenable groups with finite Hirsch length, namely wreath products of finite groups by finitely generated virtually nilpotent groups. As noted in the introduction, this not only includes many non-residually finite examples, but also non-solvable ones. 

Our strategy is built on the isomorphism $C^*(K \wr H) \cong A^{\otimes H} \rtimes H$, which leads us to study shift actions associated to multi-matrix algebras, where we use a theorem of Hirshberg and the second named author \cite{HirshbergWuActions} to reduce the problem to shift actions associated to matrix (or more generally UHF) algebras. This special case falls into the scope of a theorem of Sato \cite{Sato19} on $\mathcal{Z}$-stability, which implies finite nuclear dimension for (direct sums of) simple separable unital nuclear C*-algebras \cite{CastillejosEvingtonTikuisisWhiteWinter21}. At this last step, we need to make use of the following general fact, which is essentially \cite[Proposition 3.3]{Rieffel1980Actions}.  

\begin{lemma} \label{lemma:multi-simple-cp}
	Let $(G,A,\alpha,\omega)$ be a twisted dynamical system where $G$ is a finite group and $A$ is a finite direct sum of simple C*-algebras. Then the twisted crossed product $A \rtimes_{\alpha, \omega} G$ is also a finite direct sum of simple C*-algebras. 
\end{lemma}

\begin{proof}
	The untwisted case is immediate from \cite[Proposition 3.3]{Rieffel1980Actions}. For the general twisted case, we simply apply the Packer-Raeburn stabilization trick (see \Cref{sec:stabilization}), since a C*-algebra is a finite direct sum of simple ones if and only if its stabilization is. 
\end{proof}

In the following theorem, we write $\growthdegree(G)$ for the degree of the polynomial growth of a finitely generated virtually nilpotent group $G$. 

\begin{theorem} \label{thm:dimnuc-wr-finite}
	Let $A$ be a finite-dimensional C*-algebra and let $H$ be a finitely generated virtually nilpotent group. Then we have $\dimnuc(A^{\otimes H} \rtimes H) \leq 2 \cdot 9^{\growthdegree(H)}$, where $H \curvearrowright A^{\otimes H}$ by left translation on the tensor factors. 
\end{theorem}

\begin{proof}
	We may assume $A \not= 0$. 
	Since $A$ is finite-dimensional, it is a multi-matrix algebra, that is, we have a decomposition $A \cong \bigoplus_{i = 1}^{m} M_{n_i}(\mathbb{C})$ for some positive integers $m$ and $n_1, \ldots, n_m$. Its center can be identified with $C(\{1, \ldots, m\})$ with fibers $A_i \cong M_{n_i}(\mathbb{C})$. 
	Writing $X = \{1, \ldots, m\}^{H}$, we see that $A^{\otimes H}$ is an $H$-$X$-C*-algebra with fibers $(A^{\otimes H})_s \cong \bigotimes_{h \in H} M_{n_{s(h)}}(\mathbb{C})$ for all $s \in X$ and with $H$ acting by left translation on the indices. 
	Since $X$ is a Cantor set, we may combine the  
	inequalities \eqref{eq:bounds} and \Cref{thm:HW} to see that
	\begin{equation*}
	\dimnuc(C^*(A^{\otimes H} \rtimes H)) + 1 \leq 3^{\growthdegree(H)}\cdot 3^{\growthdegree(H)}(0 + 1)\cdot (d_{\mathrm{stab}} + 1) \; ,
	\end{equation*}
	where 
	\[
		d_{\mathrm{stab}} := \sup \left\{ \dimnuc \left((A^{\otimes H})_s \rtimes_{} N \right) \colon s \in X, N \leq H_s \right\} 
	\; .
	\]
	It remains to show $\dimnuc \left((A^{\otimes H})_s \rtimes_{} N \right) \leq 1$ for any $s \in X$ and any $N \leq H_s$. 
	
	To this end, we note that since $N \leq H_s$, for each $i \in \{1, \ldots, m\}$, the preimage $s^{-1}(i)$ is a union of right $N$-cosets and thus $(A^{\otimes H})_s$ is $N$-equivariantly isomorphic to the UHF\footnote{Here we allow UHF algebras to be finite-dimensional, in which case they are matrix algebras.} algebra $M_{n_{s, N}}(\mathbb{C})^{\otimes N}$ with the shift action, where $n_{s, N}$ denotes the supernatural number 
	\[
		\prod_{i \in \{1, \ldots, m\}} n_{i}^{|s^{-1}(i) / N|} \; .
	\]
	As described in \Cref{sec:Polycyclicgroups}, there is a characteristic subgroup $K\trianglelefteq N $ that is torsion free and has finite index. 	
	The decomposition of $N$ into $K$-cosets yields a $K$-equivariant isomorphism $M_{n_{s, N}}(\mathbb{C})^{\otimes N} \cong \left(M_{n_{s, N}}(\mathbb{C})^{\otimes N/K}\right)^{\otimes K}$. 
	Note that the UHF algebra $M_{n_{s, N}}(\mathbb{C})^{\otimes N/K}$ falls in the class $\mathcal{C}$ as defined in \Cref{def:mathcalC}. 
	Hence by \Cref{example:shift-action}, we have $M_{n_{s, N}}(\mathbb{C})^{\otimes N} \rtimes_{} K \cong \left(M_{n_{s, N}}(\mathbb{C})^{\otimes N/K}\right)^{\otimes K} \rtimes_{\alpha} K \in \mathcal{C}$, for $K$ is torsion free. 
	Now it follows from \eqref{eq:tcpastcp} that there is an isomorphism
	\begin{equation*} 
	(M_{n_{s, N}}(\mathbb{C})^{\otimes N} \rtimes_{} K)\rtimes_{\alpha,\omega}N/K \cong M_{n_{s, N}}(\mathbb{C})^{\otimes N} \rtimes_{} N 
	\end{equation*}
	for some twisted action $(\alpha,\omega)$ of the finite group $N/K$ on $M_{n_{s, N}}(\mathbb{C})^{\otimes N} \rtimes_{} K$. 
	It thus follows from \Cref{lemma:multi-simple-cp} that $M_{n_{s, N}}(\mathbb{C})^{\otimes N} \rtimes_{} N $ is a finite direct sum of simple unital C*-algebras. 
	Now if $M_{n_{s, N}}(\mathbb{C})^{\otimes N} \rtimes_{} K$ is finite dimensional, then so is $M_{n_{s, N}}(\mathbb{C})^{\otimes N} \rtimes_{} N $, which hence has nuclear dimension zero. Otherwise by \Cref{rem:equivalentCdef}, $M_{n_{s, N}}(\mathbb{C})^{\otimes N} \rtimes_{} K$ is $\mathcal{Z}$-stable and a special case of Sato's theorem \cite[Theorem 1.1]{Sato19} verifies $\mathcal{Z}$-stability for $M_{n_{s, N}}(\mathbb{C})^{\otimes N} \rtimes_{} N$ and thus also for all of its simple direct summands, which implies these simple direct summands (and thus also $M_{n_{s, N}}(\mathbb{C})^{\otimes N} \rtimes_{} N$ itself) have nuclear dimension no more than $1$ by \cite{CastillejosEvingtonTikuisisWhiteWinter21}, as desired. 
\end{proof}

\begin{corollary} \label{cor:dimnuc-wr-finite}
	Let $K$ be a finite group and let $H$ be a finitely generated virtually nilpotent group. Then we have $\dimnuc(C^*(K \wr H)) < \infty$. 
\end{corollary}

\begin{proof}
	Since $C^*(K \wr H) \cong A^{\otimes H} \rtimes H$ where $A = C^*(K)$ is a finite-dimensional algebra, the result follows directly from \Cref{thm:dimnuc-wr-finite}. 
\end{proof}

\section{Wreath products and infinite nuclear dimension} \label{sec:wr-infinite}

In this section, we turn our attention to the ``if'' direction of \Cref{conj:dimnuc} and show a large class of wreath products with infinite Hirsch length indeed give rise to group C*-algebras with infinite nuclear dimension. Combining this with \Cref{cor:dimnuc-wr-finite} and \Cref{thm:finitend}, we will complete the proof of \Cref{thm:C}. 

\begin{example} \label{example:ZwrZ}
	A prototypical example of a wreath product elementary amenable group with infinite Hirsch length is $\mathbb{Z} \wr \mathbb{Z}$. It is not hard to see that $\dimnuc(C^*(\mathbb{Z} \wr \mathbb{Z})) = \infty$, in accordance with \Cref{conj:dimnuc}. 
	Indeed, observing that $\mathbb{Z} \wr \mathbb{Z}$ has a presentation $\left\langle a, b \mid [a , b^k a b^{-k}] = 1 \text{ for any } k \in \mathbb{Z} \right\rangle$, we have, for any positive integer $n$, a surjective group homomorphism $\mathbb{Z} \wr \mathbb{Z} \to \mathbb{Z} \wr (\mathbb{Z}/n \mathbb{Z}) \cong \left\langle a, b \mid b^n = 1 = [a , b^k a b^{-k}] \text{ for any } k \in \mathbb{Z} \right\rangle$, which induces a surjective $*$-homomorphism 
	\[
		C^*(\mathbb{Z} \wr \mathbb{Z}) \to C^*(\mathbb{Z} \wr (\mathbb{Z}/n \mathbb{Z})) \cong C^*(\mathbb{Z})^{\otimes \mathbb{Z}/n \mathbb{Z}} \rtimes (\mathbb{Z}/n \mathbb{Z}) \cong C(\mathbb{T}^n)  \rtimes (\mathbb{Z}/n \mathbb{Z}) \; , 
	\]
	where $\mathbb{Z}/n \mathbb{Z}$ acts on $\mathbb{T}^n$ by permuting the factors in the Cartesian product. 
	Let $Y$ be the invariant open subset of $\mathbb{T}^n$ consisting of points with trivial stabilizers. 
	Since $\mathbb{T}^n$ is an $n$-dimensional topological manifold, so are $Y$ and $Y / (\mathbb{Z}/n \mathbb{Z})$. 
	Since $C_0(Y)  \rtimes (\mathbb{Z}/n \mathbb{Z})$ is a separable continuous trace C*-algebra, 
	it follows from \Cref{lemma:dimnuc-basic}\eqref{lemma:dimnuc-basic::extension} and~\eqref{lemma:dimnuc-basic::continuous-trace} that 
	\[
		\dimnuc(C^*(\mathbb{Z} \wr \mathbb{Z}))  \geq \dimnuc(C^*(\mathbb{Z} \wr (\mathbb{Z}/n \mathbb{Z}))) \geq  \dimnuc(C_0(Y)  \rtimes (\mathbb{Z}/n \mathbb{Z})) = \dim (Y / (\mathbb{Z}/n \mathbb{Z})) = n \; .
	\]
	As $n$ is arbitrary, we see that $\dimnuc(C^*(\mathbb{Z} \wr \mathbb{Z}))  = \infty$.

	Using the same method, one can also show $\dimnuc(C^*(\mathbb{Z}^m \wr \mathbb{Z}^n))  = \infty$. 
	
	However, it is difficult to generalize this method of establishing infinite nuclear dimension to deal with wreath products $K \wr H$ where $K$ is nonabelian or $H$ is non-residually finite. 
	While it is clear that the above argument breaks down for non-residually finite $H$, the challenge presented by nonabelian $K$ deserves more discussion. 
	Let us consider $D_{\infty} \wr \mathbb{Z}$, where $D_\infty$ is the infinite dihedral group. Due to the commutator relations in the presentation of wreath products, the surjective group homomorphism as defined above does not map onto $D_\infty \wr (\mathbb{Z}/n \mathbb{Z})$ but only onto the finite group $D_2 \wr (\mathbb{Z}/n \mathbb{Z})$, as $D_2$, the Klein group, is the abelianization of $D_\infty$. These finite quotients are useless in establishing lower bounds for $\dimnuc(C^*(D_\infty \wr \mathbb{Z}))$. 
	
	Now instead of aiming to construct a surjective group homomorphism, we work at the level of $*$-homomorphisms, by applying the $*$-isomorphisms $C^*(D_{\infty} \wr \mathbb{Z}) \cong C^*(D_\infty)^{\otimes \mathbb{Z}} \rtimes \mathbb{Z}$ and $C^*(D_\infty) \cong C^*(\mathbb{Z} \rtimes_{\mathrm{flip}} \mathbb{Z}/2\mathbb{Z}) \cong C^*(\mathbb{Z}) \rtimes_{\mathrm{flip}} \mathbb{Z}/2\mathbb{Z} \cong C(\mathbb{T}) \rtimes_{\mathrm{flip}} \mathbb{Z}/2\mathbb{Z} \cong \{ f \in C([0, \pi], M_2) \colon f(0) \text{ and } f(\pi) \text{ are diagonal} \}$. Nevertheless, if we want to generalize the above proof for $\mathbb{Z} \wr \mathbb{Z}$ that relies on the residual finiteness of the latter $\mathbb{Z}$, we will run into a subtler problem. 
	
	To highlight the issue, let us focus on 
	a homogeneous quotient C*-algebra $C([a,b], M_2)$ of $C^*(D_\infty)$ for some $0 < a < b < \pi$ 
	(using the above identification). 
	This gives rise to a quotient algebra $C([a,b], M_2)^{\otimes \mathbb{Z}} \rtimes \mathbb{Z}$ of $C^*(D_{\infty} \wr \mathbb{Z})$. 
	Now it follows from the dimension reduction phenomenon of 
	$\mathcal{Z}$-stable C*-algebras (\cite{TikuisisWinter2014}) that $C([a,b], M_2)^{\otimes \mathbb{Z}} \cong C([a,b]^\infty) \otimes M_{2^\infty}$ 
	has nuclear dimension no more than $2$, which makes $C([a,b], M_2)^{\otimes \mathbb{Z}} \rtimes \mathbb{Z}$ unlikely to have infinite nuclear dimension and thus a bad candidate for a desired lower bound of $\dimnuc(C^*(D_\infty \wr \mathbb{Z}))$. 
\end{example}

To circumvent the difficulty brought about by the presence of UHF algebras as described above, we adopt a different strategy from \Cref{example:ZwrZ}, one where we only take up finitely many matrix tensor copies at a time. 

\begin{lemma} \label{lemma:infty-group-finite-subset}
	Let $H$ be an infinite group. Then for any natural number $n$, there is a finite subset $F$ of $H$ such that $|F| \geq n$ and $h F \not= F$ for any $h \in H \smallsetminus \{e\}$. 
\end{lemma}

\begin{proof}
	Let us fix an element $g \in H \smallsetminus \{e\}$ and discuss two cases: 
	\begin{enumerate}
		\item If the order of $g$ is infinite, then we pick $F = \{g^j \colon j = 0, 1, \ldots, n-1 \}$. Now for any $h \in H \smallsetminus \{e\}$, if $h \notin \langle g \rangle$, then $h F \cap F \subseteq (h \langle g \rangle) \cap \langle g \rangle = \varnothing$ and thus $h F \not= F$; otherwise $h = g^m$ for some $m \in \mathbb{Z} \smallsetminus \{0\}$, and thus $h F = \{g^j \colon j = m, m+1, \ldots, m+n-1 \} \not= F$. 
		
		\item If the order of $g$ is equal to an integer $m > 1$, then there are infinitely many left cosets of $\langle g \rangle$, from which we pick $n$ left cosets $g_i \langle g \rangle$ for $i = 0, 1, \ldots, n-1$ with $g_0 = e$, and define 
		\[
			F = \bigsqcup_{i = 0}^{n-1} g_i \langle g \rangle \smallsetminus \{e\} = \left(\langle g \rangle \smallsetminus \{e\}\right) \sqcup \left( \bigsqcup_{i = 1}^{n-1} g_i \langle g \rangle  \right) \; .
		\]
		Thus $|F| = n |\langle g \rangle| - 1 > n$. Now for any $h \in H \smallsetminus \{e\}$, if $h \in g_i \langle g \rangle$ for some $i \in \{0, \ldots, n-1\}$, then $F \ni h = h e \notin h F$; otherwise we have $(h \langle g \rangle) \cap F = \varnothing$ and thus $h g \notin F$, though $h g \in h F$. 
	\end{enumerate}
	In all cases, we conclude that $h F \not= F$ for any $h \in H \smallsetminus \{e\}$. 
\end{proof}

\begin{lemma}\label{lemma:embedding-ideal}
	Let $H$ be an infinite discrete group. Let $A$ be a C*-algebra with $A^+$ denoting its unitization (i.e., $A^+ / A \cong \mathbb{C}$). Then for any natural number $n$, there is an embedding of $A^{\otimes n} \otimes \mathcal{K}(\ell^2(H))$ as an ideal in a quotient of $(A^+)^{\otimes H} \rtimes_{\alpha, \mathrm{r}} H$, where the action $\alpha$ of $H$ on $(A^+)^{\otimes H}$ is by left translation on the indices. 
\end{lemma}

\begin{proof}
	For a fixed natural number $n$, let $F$ be a finite subset of $H$ satisfying the conclusion of \Cref{lemma:infty-group-finite-subset}. Consider the quotient map $\pi_F \colon (A^+)^{\otimes H} \to (A^+)^{\otimes F}$ that mods out the $h$-th tensor factor $A^+$ by $A$ for all $h \in H \smallsetminus F$ (more precisely, if we realize $(A^+)^{\otimes H}$ as the direct limit of $(A^+)^{\otimes F'}$, as $F'$ ranges over all finite subsets of $H$ containing $F$, then $\pi_F$ is induced from the compatible system of quotient maps $\pi_{F}^{F'} \colon (A^+)^{\otimes F'} \cong (A^+)^{\otimes F' \smallsetminus F} \otimes (A^+)^{\otimes F} \to (A^+/A)^{\otimes F' \smallsetminus F} \otimes (A^+)^{\otimes F} \cong (A^+)^{\otimes F}$). Define an $H$-invariant ideal $I = \bigcap_{h \in H} \alpha_h \left( \ker (\pi_F) \right) \mathrel{\triangleleft} (A^+)^{\otimes H}$. Also consider an ideal $J_0 = (A^+)^{\otimes H \smallsetminus F} \otimes A^{\otimes F} \triangleleft (A^+)^{\otimes H}$ together with  
	C*-subalgebras $B = (\mathbb{C} \cdot 1_{A^+} )^{\otimes H \smallsetminus F} \otimes A^{\otimes F}$ in $J_0$ and $J = B + I$ in $(A^+)^{\otimes H}$. 
	We make a few observations: 
	\begin{enumerate}
		\item\label{lemma:embedding-ideal::proof:injective} The map $\pi_F$ is injective on $B$, whence 
		$I \cap B \subseteq \ker (\pi_F)  \cap B = 0$, that is, the quotient map $(A^+)^{\otimes H} \to (A^+)^{\otimes H} / I$ is injective on $B$, and thus also on $\alpha_h(B)$, for any $h \in H$.  
		
		\item\label{lemma:embedding-ideal::proof:technical}  For any $h \in H \smallsetminus \{e\}$, since $h F \not= F$, there is $h' \in hF \smallsetminus F$, and thus we have 
		\begin{align*}
		\pi_F(\alpha_h(J)) =  \pi_F(\alpha_h(B)) \subseteq  \pi_F(\alpha_h(J_0))  & = \pi_F \left( (A^+)^{\otimes H \smallsetminus h F} \otimes A^{\otimes h F} \right) \\
		& \subseteq \pi_F \left( (A^+)^{\otimes H \smallsetminus \{h'\}} \otimes A^{\otimes \{h'\}} \right) = \{0\} \; .
		\end{align*}
		
		\item\label{lemma:embedding-ideal::proof:disjoint}  For any $h \in H \smallsetminus \{e\}$, we have $J \cap \alpha_h(J) = I$. Indeed, as ``$\supseteq$'' is clear, it suffices to show ``$\subseteq$'', which amounts to showing 
		$\pi_F \left(\alpha_{h'} (J \cap \alpha_h(J))\right) = \{0\}$ for any $h' \in H$, but this follows directly from \eqref{lemma:embedding-ideal::proof:technical} as either $h'$ or $h'h$ will be nontrivial.  
		
		\item\label{lemma:embedding-ideal::proof:ideal}  The C*-subalgebra $J$ is an ideal in $(A^+)^{\otimes H}$. 
		Indeed, it suffices to show $I + B = I + J_0$. 
		Since $B \subseteq J_0$, this is reduced to showing $I + B \supseteq J_0$. 
		To do this, we fix $a \in J_0$ and observe that since $\pi_F(J_0) = A^{\otimes F} = \pi_F(B)$, there is $a' \in B$ such that $a - a' \in \ker (\pi_F)$. 
		Moreover, for any $h \in H \smallsetminus \{e\}$, we have $\pi_F(\alpha_h(a - a')) \in \pi_F(\alpha_h(J_0)) = \{0\}$ by \eqref{lemma:embedding-ideal::proof:technical}. 
		This shows $a - a' \in I$ and thus $a \in I + B$, as desired. 
	\end{enumerate}
	It follows from these observations that $J/I \cong B$ and the $H$-translates of $J/I$ in $(A^+)^{\otimes H} / I$ are disjoint-away-from-$0$ (and thus orthogonal) ideals. 
	Thus we obtain an $H$-equivariant embedding from the composition
	\begin{align*}
		A^{\otimes F} \otimes c_0(H) \xrightarrow{\cong} \bigoplus_{h \in H} \alpha_h (B) \xrightarrow{\cong} \bigoplus_{h \in H} \alpha_h (J / I) \xrightarrow{\cong}  \sum_{h \in H} \alpha_h (J) / I  \mathrel{\triangleleft} (A^+)^{\otimes H} / I
	\end{align*}
	where $H$ acts on $A^{\otimes F} \otimes c_0(H)$ by left translation on the second tensor factor, and the first map is determined as follows: for any $h \in H$, it maps $\left( \bigotimes_{g \in F} a_g \right) \otimes \delta_h$ to $\alpha_h \left( \bigotimes_{g \in H} a_g \right)$ inside $\alpha_h (B)$, with $a_g := 1_A$ for any $g \in H \smallsetminus F$. 
	
Taking reduced crossed products with $H$ (which preserves injectivity and surjectivity of $*$-homomorphisms), we obtain an embedding of $\left(A^{\otimes n} \otimes c_0(H)\right)  \rtimes_{\mathrm{r}} H$ as an ideal in $\left((A^+)^{\otimes H} / I\right) \rtimes_{\mathrm{r}} H$. Since the former algebra is isomorphic to $A^{\otimes n} \otimes \mathcal{K}(\ell^2(H))$ by Green's imprimitivity (see \eqref{eq:imprimitivity}), while the latter is a quotient of $(A^+)^{\otimes H} \rtimes_{\mathrm{r}} H$ induced from the $H$-equivariant quotient map $(A^+)^{\otimes H} \to (A^+)^{\otimes H} /I$, we have completed the desired construction. 
\end{proof}

We will apply \Cref{lemma:embedding-ideal} to the case where $A$ is the \emph{augmentation ideal} $I(K)$ of a discrete group $K$. Recall that $I(K)$ is the kernel of the \emph{augmentation map} $C^*(K) \to \mathbb{C}$ induced from the quotient homomorphism from $K$ to the trivial group. Thus we have $I(K)^+ \cong C^*(K)$. 

\begin{theorem} \label{thm:dimnuc-wr-infty}
	Let $H$ be an infinite group and let $K$ be a group that has an infinite, finitely generated, virtually abelian group $Q$ as a quotient. Then we have $\dimnuc(C^* (K \wr H)) = \infty$. 
\end{theorem}

\begin{proof}
	Without loss of generality, we may assume both $H$ and $K$ are amenable, as otherwise $C^* (K \wr H)$ would be non-nuclear and have infinite nuclear dimension for a trivial reason. 
	Now since the quotient homomorphism $K \to Q$ induces one from $K \wr H$ to $Q \wr H$, which further induces a surjective $*$-homomorphism from $C^* (K \wr H)$ to $C^* (Q \wr H)$, it follows from  \Cref{lemma:dimnuc-basic}\eqref{lemma:dimnuc-basic::extension} that $\dimnuc(C^* (K \wr H)) \geq \dimnuc(C^* (Q \wr H))$, whence it suffices to show the latter is infinite. 
	Since $C^* (Q \wr H) \cong C^*(Q)^{\otimes H} \rtimes H \cong (I(Q)^+)^{\otimes H} \rtimes H$, it follows from \Cref{lemma:embedding-ideal} and the permanence properties of nuclear dimension (see \Cref{lemma:dimnuc-basic}\eqref{lemma:dimnuc-basic::extension} and~\eqref{lemma:dimnuc-basic::stabilization}) that $\dimnuc(C^* (Q \wr H)) \geq \dimnuc (I(Q)^{\otimes n})$ for any positive integer $n$. Hence it suffices to show $\dimnuc (I(Q)^{\otimes n}) \to \infty$ as $n \to \infty$. 
	
	To this end, we fix a finite index abelian subgroup $M$ of $Q$. By the classification of finitely generated abelian groups and a standard trick for finite index subgroups, by shrinking $M$ if necessary, we may assume without loss of generality that $M$ is normal in $Q$ and isomorphic to $\mathbb{Z}^d$ for some positive integer $d$. 
	The conjugation action of $Q$ on $M$ thus amounts to a group homomorphism $\varphi \colon Q \to \mathrm{Aut}(M)$. 
	Let $N$ be the kernel of $\varphi$ and write $\underline{\varphi} \colon Q/N \to \mathrm{Aut}(M)$ for the induced monomorphism. 
	Note that $N$ agrees with the centralizer of $M$ and contains $M$ as a finite-index subgroup. 
	It follows from 
	\Cref{ex:beyond-central-extension} that $C^*(Q)$ decomposes as a twisted crossed product $C^*(N) \rtimes_{\alpha, \omega} Q/N$, where  $C^*(N)$ is a twisted $Q/N$-$\widehat{M}$-C*-algebra such that $Q/N \curvearrowright \widehat{M}$ via $\underline{\varphi} \colon Q/N \to \mathrm{Aut}(M) \cong \mathrm{Aut}_{\mathrm{cont}}(\widehat{M})$, where $\mathrm{Aut}_{\mathrm{cont}}(\widehat{M})$ stands for the group of continuous automorphisms on $\widehat{M}$.

	Now we pick a free orbit $Q/N \cdot x$ of the action $Q/N \curvearrowright \mathbb{T}^d$. To do this, observe that 
	$\mathrm{Aut}_{\mathrm{cont}}(\widehat{M}) \cong \mathrm{Aut}(M) \cong \mathrm{Aut} (\mathbb{Z}^d) \cong GL(d, \mathbb{Z})$ and thus $\underline{\varphi}$ embeds $Q/N$ as a subgroup of $GL(d, \mathbb{Z})$. 
	Now if we pick an arbitrary $\widetilde{x} \in \mathbb{R}^d$ such that its coordinates are $\mathbb{Q}$-linearly independent and all irrational, then $[\widetilde{x}] \in \mathbb{R}^d / \mathbb{Q}^d$ is a free point under the canonical action by $GL(d ,\mathbb{Q})$. Hence $x := [\widetilde{x}] \in \mathbb{R}^d / \mathbb{Z}^d$ is a free point under the canonical action by $GL(d ,\mathbb{Z})$ and thus in particular under the action by $Q/N$. 
	
	Since $Q/N$ is finite and $\mathbb{T}^d$ is Hausdorff, we may find an open neighborhood $U$ of $x$ in $\mathbb{T}^d \smallsetminus \{1\}$ such that $\{[q] \cdot U \colon [q] \in Q/N\}$ is a disjoint family, where $1$ denotes the trivial character of $\mathbb{Z}^d$. 
	Since by \Cref{prop:central-by-finite}, $C^*(N)$ is the C*-algebra of continuous sections on a locally trivial $C^*(N/M)$-bundle, 
	by shrinking $U$ if necessary, we may assume without loss of generality that $C_0(U) \cdot C^*(N) \cong C_0(U) \otimes C^*(N/M)$ as $U$-C*-algebras. 
	Since $\mathbb{T}^d$ is a $d$-dimensional manifold, by shrinking $U$ further if necessary, we may assume without loss of generality that $U$ is homeomorphic to $\mathbb{R}^d$. 
	Now $C_0((Q/N) \cdot U) = \bigoplus_{[q] \in Q/N} (C_0([q] \cdot U))$ is a $Q/N$-invariant ideal in $C^*(M)$, the latter being central in $C^*(Q)$. This generates the ideal $J := C_0((Q/N) \cdot U)  \cdot C^*(Q)$ in $C^*(Q)$.

	Since $1 \in \mathbb{T}^d$ is fixed by automorphisms of $\mathbb{T}^d$, we have $1 \notin (Q/N) \cdot U$. 
	Then since the composition $C^*(M) \to C^*(Q) \to \mathbb{C}$ induced from the obvious group homomorphisms agrees with the augmentation map $C^*(M) \to \mathbb{C}$, which corresponds, under the Gelfand duality, to the evaluation map at $1 \in \mathbb{T}^d$, it follows that $C_0((Q/N) \cdot U)$ and thus also $J$ are contained in $I(Q)$. 
	
	Since $\{[q] \cdot U \colon [q] \in Q/N\}$ is a disjoint family, there is a canonical $Q/N$-equivariant quotient map $(Q/N) \cdot U \to Q/N$ taking $U$ to the identity. 
	This allows us to view $C_0((Q/N) \cdot U)$ and thus also $J$ as (twisted) $Q/N$-$Q/N$-C*-algebras in the sense of \Cref{sec:G-X-algebras}. 
	It follows from \eqref{eq:imprimitivity} that 
	\begin{multline*}
		J \cong \left( C_0((Q/N) \cdot U) \cdot  C^*(N) \right) \rtimes_{\alpha, \omega} Q/N \\
		\cong \left(C_0(U) \cdot  C^*(N)\right) \otimes M_{|Q/N|}(\mathbb{C})  \cong  C_0(\mathbb{R}^{d}) \otimes C^*(N/M) \otimes M_{|Q/N|}(\mathbb{C}) 
	\end{multline*}
	Hence for any positive integer $n$, $I(Q)^{\otimes n}$ contains an ideal $*$-isomorphic to $J^{\otimes n}$, which in turn has a quotient $*$-isomorphic to $C_0(\mathbb{R}^{n d}) \otimes M_{|Q/M|^n}(\mathbb{C})$ through applying the augmentation homomorphism $F \to 1$. 
	It then follows from \Cref{lemma:dimnuc-basic}\eqref{lemma:dimnuc-basic::extension} and~\eqref{lemma:dimnuc-basic::continuous-trace} that 
	\[
	\dimnuc (I(Q)^{\otimes n}) \geq \dimnuc (J^{\otimes n}) \geq \dimnuc (C_0(\mathbb{R}^{n d}) \otimes M_{|Q/M|^n}(\mathbb{C})) = n d \xrightarrow{n \to \infty} \infty \; ,
	\]
	as desired. 
\end{proof}

We are now ready to prove the last remaining theorem in the introduction. 

\begin{proof}[{Proof of \Cref{thm:C}}]
	The first equivalence between finite Hirsch length of $K \wr H$ and finiteness of either $H$ or $K$ follows from the facts that $\hirsch(K), \hirsch(H) < \infty$ as well as the following equation derived from the additivity of Hirsch lengths \eqref{eq:Hirsch-length-additive}: 
	\[
		\hirsch(K \wr H) = \hirsch(K) \cdot \#(H) + \hirsch(H) \; ,
	\]
	where $\#(H)$ denotes the counting measure of $H$ and we adopt the convention $0 \cdot \infty = 0$. 
	
	To prove the ``only if'' direction of the second equivalence, we discuss two cases: if $K$ is finite, then $C^*(K)$ is a multi-matrix algebra and $C^*(K \wr H) \cong C^*(K)^{\otimes H} \rtimes H$, which has finite nuclear dimension by \Cref{thm:dimnuc-wr-finite}; if $H$ is finite, then $K \wr H$ contains a virtually polycyclic group $K^{\oplus H}$ as a finite index subgroup and is thus itself virtually polycyclic, whence $\dimnuc(C^*(K \wr H)) < \infty$ by \Cref{thm:finitend}. 
	
	Finally, to prove the ``if'' direction of the second equivalence, it suffices, in view of \Cref{thm:dimnuc-wr-infty}, to show whenever $K$ is infinite virtually polycyclic, it has an infinite, finitely generated, virtually abelian group $Q$ as a quotient. To see this, we follow \Cref{sec:Polycyclicgroups} to find a finite chain of characteristic subgroups $1 = G_n \trianglelefteq \ldots \trianglelefteq G_1 \trianglelefteq G_0 \trianglelefteq G$ such that $G / G_0$ is finite and $G_{i-1} / G_i$ is finitely generated abelian for $i = 1, 2, \ldots, n$. Let $G_k$ be the first in this chain such that $G_{k-1} / G_k$ is infinite. Then it follows that $G / G_k$ is infinite, finitely generated, and virtually abelian, as desired. 
\end{proof}

\appendix

\section{Central extensions and lifting homomorphisms } \label{appendix:celh}

In this appendix we digress to discuss how the central extension $1 \rightarrow H_2(G) \rightarrow E \rightarrow G \rightarrow 1$ in \Cref{lem:centralextension} depends on the choice of the lifting homomorphism $\sigma \colon B_1(G)\rightarrow C_2(G) / B_2(G)$. 
For the sake of clarity, we shall often write $\pi_\sigma$ and $E_\sigma$ in the construction in \Cref{lem:centralextension} to emphasize the dependence on $\sigma$. 

\begin{proposition} \label{remark:H2-extension-uniqueness}
	For any discrete group $G$, there is a bijection between the following sets: 
	\begin{enumerate}
		\item \label{remark:H2-extension-uniqueness::class} the set of isomorphism classes of central extensions $1 \rightarrow H_2(G) \rightarrow E \rightarrow G \rightarrow 1$ as constructed in \Cref{lem:centralextension}, 
		where two central extensions $E$ and $E'$ are considered isomorphic if there is an isomorphism $\varphi \colon E \to E'$ that makes the following diagram commute
		\[
		\xymatrix{
			1 \ar[r] & H_2(G) \ar[r] \ar@{=}[d] & E \ar[r]  \ar[d]^\varphi & G \ar[r]  \ar@{=}[d] & 1 \\
			1 \ar[r] &  H_2(G) \ar[r] &  E' \ar[r] &  G \ar[r] &  1
		}
		\]
		
		\item \label{remark:H2-extension-uniqueness::ext} the abelian group $\mathrm{Ext} (H_1(G), H_2(G))$. 
	\end{enumerate}
\end{proposition}

\begin{proof}
	The central extensions constructed in \Cref{lem:centralextension} depend solely on the choice  of lifting homomorphisms $\sigma \colon B_1(G)\rightarrow C_2(G)$. 	
	Let us fix a lifting homomorphism $\sigma \colon B_1(G)\rightarrow C_2(G) / B_2(G)$. 
	Note that any other choice of the lifting homomorphism is of the form $\sigma + \delta$ for some homomorphism $\delta \colon B_1(G) \to Z_2(G) / B_2(G) = H_2(G)$. 
	Following the construction in the proof of \Cref{lem:centralextension}, we see that 
	\begin{equation} \label{remark:H2-extension-uniqueness::eq:B1-cohomologous}
	{\pi}_{\sigma + \delta} (g_1, g_2) - {\pi}_{\sigma} (g_1, g_2) = - {\delta} \partial_2 (g_1, g_2) = - {\delta} ((g_1) - (g_1 g_2) + (g_2) ) 
	\end{equation}
	for any $g_1 , g_2 \in G$. 
	
	On the other hand, given two central extensions $1 \rightarrow H_2(G) \rightarrow E_\sigma \rightarrow G \rightarrow 1$ and $1 \rightarrow H_2(G) \rightarrow E_{\sigma + \delta} \rightarrow G \rightarrow 1$, we observe that a map $\varphi \colon E_\sigma \to E_{\sigma + \delta}$ is an isomorphism that makes the following diagram commute
	\[
	\xymatrix{
		1 \ar[r] & H_2(G) \ar[r] \ar@{=}[d] & E_\sigma \ar[r]  \ar[d]^\varphi & G \ar[r]  \ar@{=}[d] & 1 \\
		1 \ar[r] &  H_2(G) \ar[r] &  E_{\sigma + \delta} \ar[r] &  G \ar[r] &  1
	}
	\]
	if and only if there is a map $\beta \colon G \to H_2(G)$ 
	(or equivalently, a homomorphism $C_1(G) \to H_2(G)$) 
	such that 
	$\varphi(x,g) = (x + \beta(g) , g)$ and 
	\begin{equation} \label{remark:H2-extension-uniqueness::eq:C1-cohomologous}
	{\pi}_{\sigma + \delta} (g_1, g_2) - {\pi}_{\sigma} (g_1, g_2) = - \beta(g_1) + \beta(g_1 g_2) - \beta(g_2) 
	= - \beta \partial_2 (g_1, g_2) 
	\end{equation}
	for any $x \in H_2(G)$ and $g, g_1 ,g_2 \in G$.

	Comparing \eqref{remark:H2-extension-uniqueness::eq:B1-cohomologous} with  \eqref{remark:H2-extension-uniqueness::eq:C1-cohomologous}, we see that the set of isomorphism classes in \eqref{remark:H2-extension-uniqueness::class} 
	is in one-to-one correspondence with the abelian group 
	\[
	\mathrm{Hom} (B_1(G) , H_2(G)) / \iota^* (\mathrm{Hom} (C_1(G) , H_2(G))) \; ,
	\]
	where $\iota^* \colon \mathrm{Hom} (C_1(G) , H_2(G)) \to \mathrm{Hom} (B_1(G) , H_2(G))$ is the restriction homomorphism induced by the canonical embedding $\iota \colon B_1(G) \to C_1(G)$. 
	This quotient abelian group is in turn isomorphic to the abelian group $\mathrm{Ext} (H_1(G), H_2(G))$, thanks to the long exact sequence 
	\[
	0 \to \mathrm{Hom} (H_1(G) , H_2(G)) \to \mathrm{Hom} (C_1(G) , H_2(G)) \xrightarrow{\iota^*} \mathrm{Hom} (B_1(G) , H_2(G))  \to \mathrm{Ext} (H_1(G) , H_2(G)) \to 
	0
	\]
	induced from the short exact sequence $0 \to B_1(G) \xrightarrow{\iota} C_1(G) \to H_1(G) \to 0$ 
	and then truncated as the freeness of $C_1(G)$ implies that in the original long exact sequence, the term to the right of $\mathrm{Ext} (H_1(G) , H_2(G))$, namely $\mathrm{Ext} (C_1(G) , H_2(G))$, vanishes. 
\end{proof}

\begin{example} \label{ex:Klein}
	Let us illustrate \Cref{lem:centralextension} and \Cref{remark:H2-extension-uniqueness} on one of the simplest nontrivial examples, where we take $G$ to be the Klein four-group $\mathbb{Z}/2 \oplus \mathbb{Z}/2 = \langle a, b \mid a^2 = b^2 = ab = ba = 1 \rangle$. Then we have $H_1(G) \cong G$ and $H_2(G) \cong \mathbb{Z}/2$ (generated by the 2-cocycle $\omega$ with $\omega (a^i b^j , a^k b^l) = (-1)^{ij - kl}$ for any $i, j , k , l \in \mathbb{Z}$), whence $\mathrm{Ext} (H_1(G) , H_2(G)) \cong G$. 
	
	On the other hand, since $C^*(G, \omega) \cong M_2(\mathbb{C})$ is noncommutative, any extension in \Cref{lem:centralextension} must produce a nonabelian $E$. 
	Since $|E| = |G| \cdot |H_2(G)| = 8$, it is only possible that $E$ is isomorphic to either the quaternion group $Q_8 = \langle \mathbf{i},\mathbf{j} \mid \mathbf{i}^4 = \mathbf{j}^4 = 1, \mathbf{i}^2 = \mathbf{j}^2 , \mathbf{i} \mathbf{j} = \mathbf{j}^{-1} \mathbf{i} \rangle$ or the dihedral group $D_4 = \langle s,t \mid s^2 = t^4 = 1, st = t^{-1} s \rangle$. Note that $Z(Q_8) = \langle \mathbf{i}^2 \rangle \cong \mathbb{Z} / 2 \cong \langle t^2 \rangle = Z(D_4)$. 
	Up to isomorphism (of central extensions), we find four possible central extensions: 
	\begin{enumerate}
		\item $1 \to Z(Q_8) \to Q_8 \to G \to 1$ with $[\mathbf{i}] = a$ and $[\mathbf{j}] = b$; 
		\item $1 \to Z(D_4) \to D_4 \to G \to 1$ with $[s] = a$ and $[t] = b$; 
		\item $1 \to Z(D_4) \to D_4 \to G \to 1$ with $[s] = b$ and $[t] = ab$; 
		\item $1 \to Z(D_4) \to D_4 \to G \to 1$ with $[s] = ab$ and $[t] = a$.  
	\end{enumerate}
	Note that these four classes can be easily distinguished by tracking which of the three nontrivial elements in $G$ get lifted to elements of order $4$. 
	Finally, one may directly check that all of these central extensions do arise from the procedure in \Cref{lem:centralextension}. 
\end{example}

Note that in \Cref{ex:Klein}, no matter which central extension we choose, we always have $C^*(E) \cong C^*(G) \oplus C^*(G, \omega)$. 
This motivates us to 
investigate a weaker problem than the one in \Cref{remark:H2-extension-uniqueness}, namely how the C*-algebra $C^*(E_\sigma)$ depends on $\sigma$ as an $\widehat{H_2(G)}$-C*-algebra with a fixed 
$G$-grading on each fiber. 

Let us explain some notations. Given a locally compact Hausdorff space $X$ and two $X$-C*-algebras $A$ and $B$ (see \Cref{sec:G-X-algebras} with $G$ taken to be the trivial group), a $*$-homomorphism (respectively, $*$-isomorphism) $\psi \colon A \to B$ is an $X$-$*$-homomorphism  (respectively, $X$-$*$-isomorphism) if it intertwines the $*$-homomorphisms $\iota_A \colon C_0(X) \to ZM(A)$ and $\iota_B \colon C_0(X) \to ZM(B)$; 
in this case, $\psi$ induces a $*$-homomorphism $\psi_x \colon A_x \to B_x$ for any $x \in X$ with $\pi_x \circ \psi  = \psi_x \circ \pi_x$, where $\pi_x$ denotes the quotient maps $A \to A_x$ and $B \to B_x$. 

Moreover, given any discrete abelian groups $H$, we write $\mathrm{T}H$ for the torsion subgroup of $H$ and $\mathrm{F}H = H / \mathrm{T}H$ for the torsion-free quotient of $H$.

\begin{proposition} \label{remark:H2-extension-uniqueness-weak}
	For any discrete group $G$, there is a bijection between the following sets: 
	\begin{enumerate}
		\item the set of equivalence classes of lifting homomorphisms $\sigma \colon B_1(G)\rightarrow C_2(G)$, where two lifting homomorphisms $\sigma$ and $\sigma'$ are considered equivalent if there is an $\widehat{H_2(G)}$-$*$-isomorphism $\psi \colon C^*(E_\sigma) \to C^*(E_{\sigma'})$ such that for each $\theta \in \widehat{H_2(G)}$ and for any $g \in G$, $\psi_\theta \colon C^*(E_\sigma)_\theta \to C^*(E_{\sigma'})_\theta$ maps $(u_{(0, g)})_\theta$ in the domain to a scalar multiple of $(u_{(0, g)})_\theta$ in the range; 
		
		\item the abelian group $\mathrm{Ext} (H_1(G), \mathrm{F} H_2(G))$. 
	\end{enumerate}
\end{proposition}

We shall need the following lemma in the proof of this result.

\begin{lemma} \label{lemma:UC-decomposition}
	Let $H$ be a discrete abelian group $H$. Then there is a commutative diagram of abelian groups
	\begin{equation*}
	\xymatrix{
		0 \ar[r] & \mathrm{T}H \ar[r] \ar@{_(->}[d]^{\eta} & H \ar[r]^{\varphi} \ar@{_(->}[d]^{\eta} & \mathrm{F}H \ar[r] \ar@{_(->}[d]^{\underline{\eta}} & 0 \\
		0 \ar[r] & U_0(C(\widehat{H})) \ar[r] \ar[d]^{\varphi^{**}} & U(C(\widehat{H})) \ar[r]^{\lambda} \ar[d]^{\varphi^{**}} & \pi^1(\widehat{H}) \ar[r]  \ar[d]^{\varphi^{**}} & 0 \\
		0 \ar[r] & U_0(C(\widehat{\mathrm{F} H})) \ar[r] & U(C(\widehat{\mathrm{F} H})) \ar[r]^{\lambda} & \pi^1(\widehat{\mathrm{F} H}) \ar[r]  & 0 
	}
	\end{equation*}
	where 
	$\varphi \colon H \to \mathrm{F} H$ denotes the quotient map, 
	$U(C(\widehat{H}))$ is the unitary group of $C(\widehat{H})$ (that is, the set of $\mathbb{T}$-valued continuous functions on $\widehat{H}$), 
	$U_0(C(\widehat{H}))$ is the connected component of the identity in $U(C(\widehat{H}))$, 
	$\pi^1(\widehat{H})$ 
	(respectively, $\pi^1 (\widehat{\mathrm{F} H} )$) 
	is the group of homotopy classes of continuous maps from $\widehat{H}$ (respectively, $\widehat{\mathrm{F} H}$) to $\mathbb{T}$ (with pointwise multiplication giving a group structure), 
	$\eta$ embeds $H$ as $\mathrm{Hom} (\widehat{H} , \mathbb{T})$ inside $U(C(\widehat{H}))$, 
	and 
	we use $\varphi^{**}$ to denote the homomorphisms contravariantly induced from $\varphi^* \colon \widehat{\mathrm{F} H} \to \widehat{H}$, which itself is contravariantly induced from $\varphi$. 
	Moreover, the following hold true: 
	\begin{enumerate}
		\item \label{lemma:UC-decomposition::exact} The rows are exact. 
		\item \label{lemma:UC-decomposition::divisible} The groups $U_0(C(\widehat{H}))$ and $U_0(C(\widehat{\mathrm{F} H}))$ are divisible. 
		\item \label{lemma:UC-decomposition::split_injective} The composition $\varphi^{**} \circ \underline{\eta} \colon \mathrm{F}H \to \pi^1 \left(\widehat{\mathrm{F} H} \right)$ is an isomorphism; thus $\underline{\eta}$ is split injective. 
	\end{enumerate}
\end{lemma}

\begin{proof}
	We first observe that $U(C(\widehat{H}))$ is locally path-connected. It follows that $U_0(C(\widehat{H}))$ agrees with the path-connected component of the identity in $U(C(\widehat{H}))$ and thus is the kernel of the quotient map $U(C(\widehat{H})) \to \pi^1(\widehat{H})$. This shows the exactness of the middle row. Replacing $H$ by $\mathrm{F} H$, we obtain exactness of the bottom row. Since the exactness of the top row is by definition, we have \eqref{lemma:UC-decomposition::exact} holds. 
	
	Then observe that the exponentiation map is an epimorphism from the additive group $C(\widehat{H}, \mathbb{R})$ to $U_0(C(\widehat{H}))$, whence the divisibility of the latter follows from that of the former. Replacing $H$ by $\mathrm{F} H$, we see that $U_0(C(\widehat{\mathrm{F} H}))$ is also divisible, thus proving \eqref{lemma:UC-decomposition::divisible}. 
	
	Next we note that $\eta$ takes any $h \in \mathrm{T} H$ to a $\mathbb{T}$-valued continuous function with a finite range (namely a set of roots of unity), and thus a path of the form $\left( \exp \circ  (t \cdot \log) \circ \eta(h) \right)_{t \in [0,1]}$ connects $\eta(h)$ to the identity, whence $\eta(h) \in U_0(C(\widehat{H}))$. It then follows that $\underline{\eta}$ is well-defined. 
	
	It remains to show \eqref{lemma:UC-decomposition::split_injective}. 
	To do this, we first observe that if we construct $\eta' \colon \mathrm{F}H \to U(C(\widehat{\mathrm{F} H}))$ just as $\eta$ but with $\mathrm{F}H$ replacing $H$ and then write $\mathring{\eta}$ for $\lambda \circ \eta' \colon \mathrm{F}H \to \pi^1 \left(\widehat{\mathrm{F} H} \right)$, then we immediately have $\eta' \circ \varphi = \varphi^{**} \circ \eta \colon H \to U(C(\widehat{\mathrm{F} H}))$ and thus $\mathring{\eta} \circ \varphi = \lambda \circ \eta' \circ \varphi = \lambda \circ \varphi^{**} \circ \eta = \varphi^{**} \circ \underline{\eta} \circ \varphi$ by the commutativity of the diagram in the statement of the lemma. It follows that $\mathring{\eta}$ agrees with the composition $\varphi^{**} \circ \underline{\eta} \colon \mathrm{F}H \to \pi^1 \left(\widehat{\mathrm{F} H} \right)$. 
	
	Hence it suffices to show $\mathring{\eta} \colon \mathrm{F}H \to \pi^1 \left(\widehat{\mathrm{F} H} \right)$ is an isomorphism. 
	To this end, we observe that since the construction of $\mathring{\eta}$ only depends on the torsion-free quotient of $H$, we may henceforth assume without loss of generality that $H$ is torsion free. 
	Now in the special case where $H$ is finitely generated and torsion-free, we have $H \cong \mathbb{Z}^d$ and $\widehat{H} \cong \mathbb{T}^d$. A direct computation shows that $\mathring{\eta}$ takes the standard generators of $\mathbb{Z}^d$ to those of $\pi^1(\mathbb{T}^d) \cong H^1(\mathbb{T}^d) \cong \mathbb{Z}^d$, and is thus an isomorphism. 
	
	To deal with the general case, it is convenient to use a bit of category theory. Observe that taking Pontryagin duals forms a contravariant equivalence between the category $\mathbf{DTfAb}$ of discrete torsion-free abelian groups (and group homomorphisms) and the category $\mathbf{CCAb}$ of compact connected abelian groups (and continuous group homomorphisms), and thus, in particular, it takes direct limits to inverse limits and takes monomorphisms (namely, injective homomorphisms) to epimorphisms (namely, surjective continuous homomorphisms). 
	Meanwhile, $\pi^1$ also forms a contravariant functor from $\mathbf{CCAb}$ to the category $\mathbf{DAb}$ of discrete abelian groups (and group homomorphisms). 
	Moreover, one verifies directly that $\mathring{\eta}$ forms a natural transformation from the inclusion functor $\mathbf{DTfAb} \to \mathbf{DAb}$ to the composition of $\pi^1$ with the functor of taking Pontryagin duals. 
	
	Now observe that any torsion-free $H$ is the union of all of its finitely generated (necessarily abelian and torsion-free) subgroups. 
	Note that once we direct the set $I$ of all finitely generated subgroups by inclusion, this union amounts to the inductive limit in $\mathbf{DTfAb}$ and also in $\mathbf{DAb}$. 
	Hence, given the naturality of $\mathring{\eta}$ and what we have proved in the case of finite generated torsion-free groups, it suffices to show $\pi^1(\widehat{H})$ is isomorphic to the inductive limit of $\varinjlim_{K \in I} \pi^1(\widehat{K})$. Note that by universality, we already have a canonical homomorphism $\Phi \colon \varinjlim_{K \in I} \pi^1(\widehat{K}) \to  \pi^1(\widehat{H})$. It suffices to show $\Phi$ is an isomorphism. 
	
	As observed above, $\widehat{H}$ is homeomorphically isomorphic to the inverse limit $\varprojlim_{K \in I} \widehat{K}$, and for each $K \in I$, the canonical homomorphism $\widehat{H} \to \widehat{K}$ is surjective. 
	Now to show the surjectivity of $\Phi$, we note that given any continuous function $f \colon \widehat{H} \to \mathbb{T}$, by the construction of the topology on $\widehat{H} \cong \varprojlim_{K \in I} \widehat{K}$ and a compactness argument, there exists some $K_0 \in I$ such that when restricted on any fiber of the epimorphism $\widehat{H} \to \widehat{K_0}$, the range of $f$ is contained in a semicircle, which allows us to define 
	a continuous function $F \colon \widehat{K_0} \times [0,1] \to \mathbb{T}$ by 
	\[
		F(\chi, t) = \exp \left( t \, \log f(\chi) + (1-t)  \int_{L} \log f (\chi + \rho) \, d \mu(\rho) \right) \quad \text{ for any } \chi \in \widehat{H} \text{ and } t \in [0,1] \; ,
	\]
	where $L = \ker (\widehat{H} \to \widehat{K_0})$, $\mu$ is the normalized Haar measure on $L$, and $\log$ is a branch of the logarithmic function on a semicircle containing $\{ f(\chi + \rho) \colon \rho \in L \}$ 
	(it does not matter which branch we pick). 
	Note that $f = F(-, 1)$ while $F(-,0)$ factors through $\widehat{K_0}$. This proves surjectivity. Applying the same technique to continuous functions $\widehat{H} \times [0,1] \to \mathbb{T}$ then proves injectivity. 
\end{proof}

\begin{proof}[Proof of \Cref{remark:H2-extension-uniqueness-weak}]
	Fix a lifting homomorphism $\sigma \colon B_1(G) \to C_2(G) / B_2(G)$. As noted before, the set of all lifting homomorphisms is nothing but $\{\sigma + \delta \colon \delta \in \mathrm{Hom}(B_1(G) \to H_2(G))\}$. Fix such a $\delta$. 
	Observe that equivalence between $\sigma$ and ${\sigma + \delta}$ amounts to the existence of a map $\beta \colon G \to U(C(\widehat{H_2(G)}))$, or equivalently, a homomorphism $\beta \colon C_1(G) \to U(C(\widehat{H_2(G)}))$, 
	such that $\psi (u_{(0, g)}) = \beta(g) \cdot u_{(0, g)}$ for any $g \in G$ and 
	\begin{equation} \label{remark:H2-extension-uniqueness-weak::eq:U-cohomologous}
	\theta \left( {\pi}_{\sigma + \delta} (g_1, g_2) - {\pi}_{\sigma} (g_1, g_2) \right)= \beta(g_1)_\theta ^* \, \beta(g_1 g_2)_\theta \, \beta(g_2)_\theta ^* 
	= \beta ( \partial_2 (g_1, g_2) )_\theta ^*
	\end{equation}
	for any $g_1, g_2 \in G $ and any $\theta \in \widehat{H_2(G)}$. 
	
	Comparing \eqref{remark:H2-extension-uniqueness::eq:B1-cohomologous} with \eqref{remark:H2-extension-uniqueness-weak::eq:U-cohomologous}, we see that $\sigma$ and ${\sigma + \delta}$ are equivalent if and only if there exists a homomorphism $\beta \colon C_1(G) \to U(C(\widehat{H_2(G)}))$ that fits into the following commuting diagram 
	\begin{equation*}
	\xymatrix{
		B_1(G) \ar[r]^{{\delta}} \ar@{_(->}[d]^{\iota} & H_2(G) \ar@{_(->}[d]^{\eta} \\
		C_1(G) \ar[r]^{\beta} & U(C(\widehat{H_2(G)}))
	}
	\end{equation*}
	where $\eta$ is as in \Cref{lemma:UC-decomposition}.  
	It follows that the set of equivalent classes of lifting homomorphisms is in one-to-one correspondence to the quotient abelian group 
	\[
	\eta_* \left( \mathrm{Hom} ( B_1(G) , H_2(G) ) \right) / \left( \eta_* \left( \mathrm{Hom} ( B_1(G) , H_2(G) ) \right) \cap \iota^* \left( \mathrm{Hom} ( C_1(G) , U(C(\widehat{H_2(G)})) )  \right) \right)
	\]
	which we see is in turn isomorphic to the group 
	\begin{equation} \label{remark:H2-extension-uniqueness-weak::eq:answer}
	\mathrm{Im} \left( \eta_* \colon \mathrm{Ext} (H_1(G), H_2(G)) \to \mathrm{Ext} (H_1(G), U(C(\widehat{H_2(G)})))  \right)
	\end{equation}
	by chasing around the following commutative diagram with exact rows 
	\[
	\xymatrix{
		\mathrm{Hom} ( C_1(G) , H_2(G) ) \ar[r]^{\iota^*} \ar[d]^{\eta_*}  & \mathrm{Hom} ( B_1(G) , H_2(G) ) \ar@{->>}[r]  \ar[d]^{\eta_*} & \mathrm{Ext} ( H_1(G) , H_2(G) ) \ar[r] \ar[d]^{\eta_*}  &  0 \\
		\mathrm{Hom} ( C_1(G) , U(C(\widehat{H_2(G)})) ) \ar[r]^{\iota^*}  & \mathrm{Hom} ( B_1(G) , U(C(\widehat{H_2(G)})) ) \ar@{->>}[r]  & \mathrm{Ext} ( H_1(G) , U(C(\widehat{H_2(G)})) ) \ar[r]  & 0
	}
	\]
	where we used the fact that $C_1(G)$ is free abelian to derive that the rightmost terms, which a priori are Ext groups, are in fact $0$. 
	
	Now, the commutative diagram in \Cref{lemma:UC-decomposition} gives rise to the following commutative diagram 
	\[
	\xymatrix{
		\mathrm{Ext} ( H_1(G) , \mathrm{T} H_2(G) ) \ar[r]^{\iota^*}  & \mathrm{Ext} ( H_1(G) , H_2(G) ) \ar@{->>}[r]  \ar[d]^{\eta_*} & \mathrm{Ext} ( H_1(G) , \mathrm{F} H_2(G) ) \ar[r] \ar@{^(->}[d]^{\underline{\eta}_*}  &  0 \\
		0 \ar[r]^{\iota^*}  & \mathrm{Ext} ( H_1(G) , U(C(\widehat{H_2(G)})) ) \ar[r]^{\cong}  & \mathrm{Ext} ( H_1(G) , \pi^1(\widehat{H_2(G)}) ) \ar[r]  & 0 
	}
	\]
	with all rows exact (each being a part of a long exact sequence involving the Ext groups), where we used the divisibility of $U_0(C(\widehat{H}))$ to replace $\mathrm{Ext} ( H_1(G) , U_0(C(\widehat{H_2(G)})) )$ by $0$ in the lower leftmost term. 
	It follows from the split injectivity of $\underline{\eta}$ in \Cref{lemma:UC-decomposition} that $\underline{\eta}_*$ in the above diagram is also (split) injective. 
	Hence chasing around the above diagram, we see that the group in \eqref{remark:H2-extension-uniqueness-weak::eq:answer} is isomorphic to $\mathrm{Ext} ( H_1(G) , \mathrm{F} H_2(G) )$. 
\end{proof}

\bibliographystyle{alpha}

\bibliography{NDpolycyclic}

\end{document}